\documentclass[ijoc,nonblindrev]{informs3-neutral} 

\OneAndAHalfSpacedXII 

\usepackage{natbib}
 \bibpunct[, ]{(}{)}{,}{a}{}{,}%
 \def\bibfont{\small}%
 \def\BIBand{and}%

\usepackage[hidelinks]{hyperref}       
\usepackage{url}            
\usepackage{booktabs}       
\usepackage{amsfonts}       
\usepackage{nicefrac}       
\usepackage{microtype}      
\usepackage{lipsum}
\usepackage{graphicx}
\usepackage{indentfirst}

\usepackage{amssymb}
\usepackage{mathtools}
\usepackage{bbm}
\usepackage{centernot, cancel}
\usepackage[overload]{empheq}
\usepackage{cleveref} 

\usepackage[boxed,linesnumbered,noend]{algorithm2e}
\DontPrintSemicolon
\SetAlCapSkip{0.5em}
\SetKwComment{tcp}{\hfill $\triangleright$\ }{}
\SetKwFor{While}{while}{}{end}%
\SetKwFor{For}{for}{}{end}%

\usepackage{booktabs}
\usepackage{graphicx}
\usepackage[table,xcdraw]{xcolor}
\usepackage{float}
\usepackage{caption}
\usepackage{subcaption}
\usepackage{bm}

\TheoremsNumberedThrough     

\EquationsNumberedThrough    

\DeclarePairedDelimiter{\ceil}{\lceil}{\rceil}
\DeclarePairedDelimiter{\floor}{\lfloor}{\rfloor}

\begin{document}


\RUNAUTHOR{Ferreira, Pessoa, and Vidal}

\RUNTITLE{Influence Optimization in Networks: New Formulations and Valid Inequalities}

\TITLE{Influence Optimization in Networks:\linebreak New Formulations and Valid Inequalities}

\ARTICLEAUTHORS{%
\AUTHOR{Vinicius Ferreira}
\AFF{Departamento de Inform\'{a}tica, Pontif\'{i}cia Universidade Cat\'{o}lica do Rio de Janeiro, 
\EMAIL{vsouza@inf.puc-rio.br}}
\AUTHOR{Artur Pessoa}
\AFF{N\'{u}cleo de Logistica Integrada e Sistemas LOGIS, Universidade Federal Fluminense, 
\EMAIL{artur@producao.uff.br}}
\AUTHOR{Thibaut Vidal}
\AFF{CIRRELT \& SCALE-AI Chair in Data-Driven Supply Chains, Department of Mathematics and Industrial Engineering, Polytechnique Montréal, Canada\\
Departamento de Inform\'{a}tica, Pontif\'{i}cia Universidade Cat\'{o}lica do Rio de Janeiro
\EMAIL{thibaut.vidal@cirrelt.ca}}
} 

\ABSTRACT{%
Influence propagation has been the subject of extensive study due to its important role in social networks, epidemiology, and many other areas. Understanding propagation mechanisms is critical to control the spread of fake news or epidemics. In this work, we study the problem of detecting the smallest group of users whose conversion achieves, through propagation, a certain influence level over the network, therefore giving valuable information on the propagation behavior in this network. We develop mixed integer programming algorithms to solve this problem. The core of our algorithm is based on new valid inequalities, cutting planes, and separation algorithms embedded into a branch-and-cut algorithm. We additionally introduce a compact formulation relying on fewer variables. Through extensive computational experiments, we observe that the proposed methods can optimally solve many previously-open benchmark instances, and otherwise achieve small optimality gaps. These experiments also provide various insights into the benefits of different mathematical formulations.
}%



\maketitle

%


\section{Introduction}
\label{sec:Intro}

Social networks shape how we communicate information, news, rumors, and opinions. Due to their dramatic impacts on society and economic importance, the mechanisms driving the spread of influence on social networks have been the focus of extensive study, motivated by applications in viral marketing  \citep{brown1987social,goldenberg2001talk,chen2010scalable,demaine2014influence, gunnecc2020least} and influence maximization  \citep{weng2010twitterrank, bakshy2011everyone}. Proper understanding of propagation mechanisms is also becoming crucial for designing efficient control measures against fake news, harmful information, and even epidemics \citep{kimura2009blocking, dreyer2009irreversible, zhang2016misinformation, zhang2019proactive}.

The recent survey of \citet{Banerjee2020} exemplifies the propagation process on the example of a commercial establishment seeking to attract new customers. A method to foster the adoption of a new product consists in distributing free samples and discounts to a group of influential individuals. If these individuals like the product, they will likely share the information with their connections, initiating a cascading process which has the potential of reaching a significant fraction of the users. As the number of free samples or financial incentives is limited, the effectiveness of this cascading process largely depends on a strategic choice of initial individuals. However, since any combination of the $n$ initial influencers represents a possible strategy, the number of solutions to this problem grows as~$2^n$, giving rise to a combinatorial problem presenting significant computational challenges.

These methodological challenges have stimulated several studies with an optimization viewpoint. Most notably, \citet{fischetti2018least} recently studied a general setting called the Generalized Least Cost Influence Problem (GLCIP). In this model, the individuals are represented by vertices in a network and connected by (possibly asymmetric) arcs representing influence relations. The goal of the GLCIP is to identify a minimum group of individuals that must be initially converted through incentives to achieve a predefined target level of influence across the network. This model allows partial incentives and permits modeling diminishing, proportional, or increasing marginal influences from multiple sources through linear or nonlinear activation functions. Its generality permits it to represent a vast array of applications but also leads to significant computational hurdles. To solve it, \citet{fischetti2018least} proposed solution algorithms based on two mixed integer programming (MIP) formulations: a formulation based on arc-flow variables, and a set-covering formulation. The study also proposed new classes of valid inequalities that permit to improve the linear relaxation. Together, these methods could identify optimal solutions for cases where the degree of the vertices (i.e., the number of connections) remains relatively small. Still, vast improvements remain needed to solve cases where individuals are more numerous or have a larger number of connections.

To contribute towards circumventing these computational issues, we propose new mathematical formulations and valid inequalities for the GLCIP. We first show how to extend the arc-flow model presented in \citet{fischetti2018least} to deal with nonlinear activation functions not initially covered by it. Then, we introduce new cutting planes and separation algorithms and use them to devise an efficient branch-and-cut algorithm. Additionally, we introduce a new compact formulation that uses only one set of variables to represent the incentive paid to each individual (possibly zero), adapting the valid inequalities and separation algorithms introduced to this case. We evaluate the method on a set of synthetic instances provided by \citet{fischetti2018least}. Our results show that the algorithm can reach the optimal value in more than 150 open cases. For the other cases, it achieves smaller optimality gaps than previously-developed approaches.

Therefore, the main contributions of this paper are:
\begin{itemize}
    \item A generalization of the polynomial formulation presented in \citet{fischetti2018least} to deal with cases with nonlinear activation functions;
    \item New valid inequalities inspired by the generalized propagation constraints proposed in \citet{fischetti2018least}, with an approach to express these inequalities within a simpler formulation based on arcs, instead of a formulation based on sets;
    \item A new compact formulation with only one set of variables to represent the incentives paid to individuals, as well as adaptations of the valid inequalities and separation algorithms to this case;
    \item Extensive computational experiments performed on the set of instances presented in \citet{fischetti2018least} which demonstrate that (i) the proposed methods reduce the gap between the lower and upper bounds of the value of the best solution, and (ii) these algorithms are capable of optimally solving a large number of additional open instances.
\end{itemize}

The remainder of this paper is structured as follows: Section~\ref{sec:Problem} formally defines the problem and reviews related works. Section~\ref{sec:Metho} describes our new formulations and valid inequalities. Section~\ref{sec:Experiments} presents our numerical analyses, and Section~\ref{sec:Conclusions} finally concludes.

\section{Related Works and Models}
\label{sec:Problem}

\citet{domingos2001mining} and \citet{richardson2002mining} were among the first to study influence propagation from an optimization viewpoint. The authors described what would later be known as the influence maximization problem in the context of viral marketing. In their formulation, a small set of influential users must be targeted initially to trigger chain adoption. Influence is modeled within a probabilistic framework using Markov random fields. 

Later on, \citet{kempe2003maximizing} revisited the influence maximization problem. The authors adopted the linear threshold diffusion model of \cite{granovetter1978threshold} to represent the propagation mechanism step by step. Given a directed graph, the goal of the problem is to select $k$ seed individuals to maximize the total propagation in the network according to the diffusion model. However, using such a model for influence maximization leads to significant computational hurdles. \citet{kempe2003maximizing} proved that computing the optimal subset of initial seeds is NP-Hard and proposed approximation algorithms.

Several studies followed this research line and extended the original model. In particular, heterogeneous costs for activating different nodes in the network are considered in the Target Set Selection Problem (TSSP -- \citealt{chen2009approximability}) and the Weighted Target Set Selection Problem (WTSSP -- \citealt{zhang2016mathematical}). In the Least Cost Influence Problem (LCIP) of \citet{gunnecc2020least}, partial incentives can be applied on a node, allowing the combination of incentives and the influence of neighbors. Finally, the Generalized Least Cost Influence Problem (GLCIP), proposed by \citet{fischetti2018least} further generalizes the model with possible diminishing, proportional, or increasing marginal influences from multiple sources.\\

\noindent
\textbf{Generalized Least Cost Influence Problem.}
The GLCIP can be modeled using a directed graph $G = (V,A)$, in which the vertices represent the individuals in the network, and the arcs map the relationships between them. Each $(i,j) \in A$ represents an influence relation from node $i \in V$ over $j \in V$ with weight $d_{ij}$. Each node has a set of possible incentives $P_{i}$ (e.g., financial incentives such as discounts or free products) that can be offered to influence it, and each incentive has an associated cost $w_{ip} \geq 0$ for offering incentive $p \in P_{i}$ to node $i \in V$. Moreover, each node $i \in V$ of this graph has an activation threshold $h_{i} > 0$, which defines the limit that needs to be attained through a combination of influences and incentives offered to activate that node. 

The set $N_{i} = \{j \in V | (j,i) \in A\}$ denotes the set of neighbors of a node $i \in V$, that is, the nodes that directly influence node $i \in V$. Let $U \subseteq N_{i}$ be the set of active neighbors of node~$i$, $p \in P_{i}$ be the incentive $p$ offered to node $i$ and $f_{i}$ a activation function of node $i$. Then, this previously inactive node $i \in V$ becomes active (influenced) at the current step of the propagation process if and only if $f_{i}(U, p) \geq h_{i}$, i.e., if a possible incentive offered $p$ and the influence received from a subset of active neighbors $U \subseteq N_{i}$ attains its threshold $h_{i}$.

\citet{fischetti2018least} use an activation function defined as $f_{i}(U, p)~=~\left(\sum_{j \in U} {d_{ji}}\right)^ {\Gamma}~+~p$. When $\Gamma = 1$, $f_{i}$ is linear and additively separable. This means that the effect of one variable on the value of the function does not depend on the values of the other variables. This also greatly simplifies the analysis \citep{buja1989}. In contrast, values of $\Gamma> 1$ model situations of peer pressure, when the marginal influence of each incentive grows with the number of incentives. Finally, $\Gamma < 1$ characterizes a submodular function and models the situation of decreasing marginal influences \citep{fischetti2018least}. In this work, we follow the same convention and use the same activation function.
Thus, the objective of GLCIP is to determine which incentive $p_i$ should be offered to each node $i \in V$ to minimize total costs $\sum_{i \in V}{w_{i,p_{i}}}$, while simultaneously ensuring that a predetermined fraction of nodes $0 \leq \alpha \leq 1$ is active at the end.

\citet{fischetti2018least} introduced two mixed integer programming formulations for the GLCIP: a formulation using arc variables (ARC), which generalizes similar arc formulations developed for other GLCIP variants, and a set-covering (COV) formulation with an exponential number of variables. Two variations of the COV model are suggested, using different separation algorithms (heuristic and exact) for a family of cuts. Based on the numerical results reported by the authors, the ability to solve the GLCIP to optimality depends to a large extent on the degree of the vertices in the network, as the proposed formulations are more efficient for sparse graphs with small average degree.\\

\noindent
\textbf{Set covering formulation and valid inequalities for the GLCIP.}
The COV model of \citet{fischetti2018least} uses the notion of \emph{Minimal Influencing Set (MIS)}. Let $U \subseteq N_{i}$ be a (possibly empty) subset of active neighbors of node $i \in V$, and let $p_{i}^{U} = \min \left \{\tilde{p} \in P_{i} ~ | ~ f_{i} (U, \tilde{p}) \geq h_{i} \right \}$ be the smallest incentive (possibly equal to zero) required for activation with these neighbors, with an incentive cost of $w_{i}^{U}$. $U$ is a MIS if and only if there is no smaller subset $U' \subset U$ such that $f_{i} (U', p_{i}^{U}) \geq h_{i}$.
With this definition, let $\Lambda_{i}$ be the set of all MIS for each node~$i \in V$.

They define three sets of variables:
\begin{itemize}
    \item Node variables $x_{i} \in \{0, 1\}$, taking value $1$ if and only if $i \in V$ is activated;
    \item Arc variables $z_{ij} \in \{0, 1\}$, taking value $1$ if and only if arc $(i, j) \in A$ carries influence;
    \item Influencing subset variables $\lambda_{i}^{U} \in \{0, 1\}$, taking value $1$ if and only if influencing set $U \in \Lambda_{i}$ is used to activate node $i \in V$.
\end{itemize}

With this, the (COV) model can be defined as follows:
\begin{align}[left = \hspace*{-1cm}\text{(COV)}]
\min \hspace*{0.4cm} & \sum_{i \in V}{\sum_{U \in \Lambda_{i}}{w_{i}^{U} \lambda_{i}^{U}}} \label{model:cov_obj}\\
\text{s.t.} \hspace*{0.4cm} & \sum_{U \in \Lambda_{i}}{\lambda_{i}^{U}} \geq x_{i}, & \forall i \in V \label{model:cov_propagation_constraint} \\
& \sum_{\substack{
  U \in \Lambda_{j} \\
  i \in U
}}{\lambda_{j}^{U}} = z_{ij}, & \forall (i,j) \in A \label{model:cov_arcs_influence} \\
& \sum_{(i,j) \in C}{z_{ij}} \leq \sum_{i \in V(C)\backslash\{k\}}{x_{i}}, & \forall \, \text{cycle} ~ C \in A, \forall k \in V(C)  \label{model:cov_constraint_cycle_cut}\\
& z_{ij} \leq x_{i}, & \forall (i,j) \in A  ~ | ~ (j,i) \notin A \label{model:cov_constraint_linking}\\
& \sum_{i \in V}{x_{i}} \geq \ceil{\alpha|V|} \label{model:cov_constraint_coverage} &\\
& x_{i} \in \{0, 1\}, & \forall i \in V \label{model:cov_constraint_x}\\
& z_{ij} \in \{0, 1\}, & \forall (i,j) \in A \label{model:cov_constraint_z} \\
& \lambda_{i}^{U} \geq 0, & \forall i \in V, \forall U \in \Lambda_{i}. \label{model:cov_constraint_lambda}
\end{align}

In this model, Objective~(\ref{model:cov_obj}) minimizes the total cost of all incentives offered. Constraints~(\ref{model:cov_propagation_constraint}) ensure that at least one subset of influences is selected for each activated node~\mbox{$i \in V$}. Constraints~(\ref{model:cov_arcs_influence}) ensure consistency between the influence subsets and arcs.
Constraints~(\ref{model:cov_constraint_cycle_cut}) ensure that the sub-graph formed by the influence arcs is acyclic ($V(C) = \{i \in V ~ | ~ (i, j) \in C \}$ refers to the nodes of the cycle). Constraints (\ref{model:cov_constraint_linking}) ensure that the influence can only be exerted on arc $(i, j) \in A$  if the node $i \in V$ is active. Finally, Constraint~(\ref{model:cov_constraint_coverage}) ensures that a proportion of at least $\alpha$ nodes becomes active.

Note that the number of $\lambda_{i}^{U}$ variables representing the possible MIS grows exponentially with $|V|$. Note that the exact methods proposed in  \citet{fischetti2018least} that use (COV) rely on the enumeration of all such variables, what explains its worse performance for graphs with large degrees. Moreover, there is also an exponential number of constraints~(\ref{model:cov_constraint_cycle_cut}), so a dedicated separation algorithm is needed to add these constraints to the model dynamically. This algorithm is described in Section \ref{sec:cycle_elimination_sep}.\\

\noindent
\textbf{Generalized propagation constraints.}
\citet{fischetti2018least} also proposed a family of cutting planes for the COV model, called \emph{generalized propagation constraints}, which strengthens the model (i.e., leads to better linear-relaxation bounds). These constraints dominate the cycle elimination constraints of the original formulation. These cutting planes are based on the observation that for each set of nodes $X \subseteq V$ containing an active node $k$, at least one of the following conditions must be true:
\begin{itemize}
    \item {At least one node $j \in X$ is activated through a non-empty minimal influencing set $U \in \Lambda_{j}$ that does not contain nodes from $X$, i.e., $U \cap X = \varnothing$;}
    \item {At least one node $j \in X$ receives $p_{j}^{\varnothing} ~ | ~ f_{j} (\varnothing, p_{j}) \geq h_{j} $, so this node $j$ is activated without receiving any influence from its neighbors $N_{j}$.}
\end{itemize}

If none of these conditions are valid, then the subgraph induced by the set of chosen arc variables induces a cycle inside $X$. Thus, the following inequalities are valid:
\begin{align}
    \sum_{j \in X}{\sum_{\substack{
        U \in \Lambda_{j}\\
        U \cap X=\varnothing
        }}{\lambda_{j}^{U} \geq x_{k}}} & & \forall k \in X, \forall X \subseteq V. \label{constraint:generalized_propagation_constraints}
\end{align}

Note that both conditions can be covered allowing $\varnothing \in \Lambda_j$. \citet{fischetti2018least} introduced two algorithms to separate Constraints~(\ref{constraint:generalized_propagation_constraints}). One is a heuristic, whereas the other is exact, leading to two variants of the COV model called C$^{+}$ and C$^{+}_{e}$, respectively.

\section{Proposed Methodology}
\label{sec:Metho}

This section introduces new mathematical programming models proposed to solve GLCIP.
We start from an arc-based formulation originally presented in \citet{fischetti2018least} and present a generalization of this arc-based formulation to deal with cases with $\Gamma \neq 1$. Next, we introduce new cutting planes and separation algorithms and incorporate them into a branch-and-cut algorithm. Finally, we propose a new compact formulation and adapt the previously developed cutting planes and separation algorithms.

\subsection{Arc Flow Formulation}
\label{sec:arc_formulation}

The arc flow formulation (ARC) emerged from earlier works on influence maximization, and it was adapted by \citet{fischetti2018least} to the GLCIP with activation function of the form $f_{i}(U, p) = \sum_{j \in U} {d_{j i}} + p$. The main characteristic of this formulation is that it guarantees a directed acyclic propagation graph (DAG). In this model, the variables $x_{i} \in \{0, 1 \}$ for $i \in V$ indicate when the node $i$ is active (influenced) or not; the arc variables $z_{ij} \in \{0, 1 \}$ for $(i, j) \in A$ indicate in which arcs the influence is exerted, that is, $z_{ij} = 1$ if the node $i$ influences the node $j$ and $0$ otherwise; and the variables $y_{ip} \in \{0, 1 \}$ for $i \in V$ and $p \in P_{i}$ refer to the incentive $p$ given to node $i$.
\begin{align}[left = \hspace*{-1cm}\text{(ARC)}]
\min \hspace*{0.4cm} & \sum_{i \in V}{\sum_{p \in P_{i}}{w_{ip} y_{ip}}} \label{model:obj}\\
\text{s.t.} \hspace*{0.4cm} & \sum_{p \in P_{i}}{py_{ip}} + \sum_{j:(j,i) \in A}{d_{ji} z_{ji}} \geq h_{i} x_{i} & \forall i \in V \label{model:constraint_propagation} \\
& \sum_{p \in P_{i}}{y_{ip}} = x_{i} & \forall i \in V \label{model:constraint_one_incentive} \\
& \sum_{(i,j) \in C}{z_{ij}} \leq \sum_{i \in V(C)\backslash\{k\}}{x_{i}} & \forall ~ \text{cycle} ~ C \in A, \forall k \in V(C)  \label{model:constraint_cycle_cut}\\
& z_{ij} \leq x_{i} & \forall (i,j) \in A ~ | (j,i) \notin A \label{model:constraint_linking}\\
& \sum_{i \in V}{x_{i}} \geq \ceil{\alpha ~ |V|} \label{model:constraint_coverage} &\\
& x_{i} \in \{0, 1\} & \forall i \in V \label{model:constraint_x}\\
& y_{ip} \in \{0, 1\} & \forall i \in V, ~ \forall p \in P_{i} \label{model:constraint_y} \\
& z_{ij} \in \{0, 1\} & \forall (i,j) \in A \label{model:constraint_z}
\end{align}

The objective function of (\ref{model:obj}) minimizes the sum of all costs related to the incentives offered. Constraints (\ref{model:constraint_propagation}) model the spread of influence between nodes. If the influences received from neighbors plus the direct incentive reach the activation threshold $h_i$, then  node~$i$ becomes active. Equations~(\ref{model:constraint_one_incentive}) ensure that each active node receives exactly one (possibly equal to zero) incentive. Constraints~(\ref{model:constraint_cycle_cut}) are used to eliminate cycles. Constraint~(\ref{model:constraint_linking}) guarantees that a node $i$ can only influence a node $j$ if the node $i$ is active. Finally, Constraint~(\ref{model:constraint_coverage}) sets up a target on the minimum fraction of nodes that need to be active.

It is worth mentioning that in \citet{fischetti2018least}, the activation function is given by $f_{i}(U, p) = \left (\sum_{j \in U} {d_{ji}} \right)^{\Gamma} + p$, and thus Constraint~(\ref{model:constraint_propagation}) only works for $\Gamma = 1$. In Section~\ref{sec:extension_gamma}, we will introduce a generalization of this constraint which is applicable for other values of~$\Gamma$.

\subsection{Separation of the Generalized Cycle Elimination Constraints}
\label{sec:cycle_elimination_sep}

As in the (COV) formulation, there is an exponential number of Constraints~(\ref{model:constraint_cycle_cut}), and a separation algorithm is needed to add them to the model iteratively. The algorithm to separate these constraints remains the same as in \citet{fischetti2018least}. It is an adaptation of the shortest path algorithm proposed in \citet{grotschel1985acyclic} to separate the classic cycle elimination cuts. Let $\bar{x}_{i}$ and $\bar{z}_{ij}$ refer to the relaxed values of the variables (\ref{model:constraint_x}) and (\ref{model:constraint_z}) respectively. Each arc $(i, j)$ of the graph receives a weight $\bar{w}_{ij} = \bar{x}_{i} - \bar{z}_{ij} \geq 0$. A shortest path between each node $k \in V$ and all of its neighbors $N_{k}$ is computed using Floyd-Warshall algorithm \citep{floyd1962algorithm, warshall1962theorem}. If the total weight of any cycle $C$ formed by the path from $k$ to a neighbor node $j \in N_{k}$ together with arc $(j, k)$ is less than $\bar{x}_{k}$, then we have found a violated cycle elimination constraint (\ref{model:constraint_cycle_cut}) which can be added to the model.

\subsection{Extension of the ARC Model}
\label{sec:extension_gamma}

In the general version of the GLCIP, the activation of a node is driven by function $f_{i}(U, p)~=~\left(\sum_{j \in N_{i}} {d_{ji}} \right)^{\Gamma} + p$. When $\Gamma = 1$, $f_{i}$ is a linear and additively separable function, this means that the effect of one variable on the value of the function does not depend on the values of the other variable, greatly simplifying the analysis \citep{buja1989}. In contrast, $\Gamma> 1$ corresponds to peer-pressure situations, whereas $\Gamma < 1$ models situations of diminishing marginal influences, which are both important and practical settings.

So far, the (ARC) formulation could only be applied to the linear case where $\Gamma = 1$. In this subsection, we modify the model to allow different values of $\Gamma$ despite the nonlinearity of the resulting activation function.

To achieve this, we start by including $\Gamma$ in Constraint (\ref{model:constraint_propagation}), leading to:
\begin{align}
    \sum_{p \in P_{i}}{p y_{ip}} + \left(\sum_{j: (j,i) \in A}{d_{ji} z_{ji}}\right)^{\Gamma} \geq h_{i} x_{i}, & & \forall i \in V. \label{new_constraint_gamma} &
\end{align}

This constraint cannot be effectively used in its current form due to its nonlinearity. However, we remark that only one incentive can be offered for each node $i \in V$. Therefore, the sum $\sum_{p \in P_{i}} {p y_{i,p}}$ will only have one incentive $p'$ such that $y_{i, p'} = 1$ and $y_{i, p} = 0, \forall p \in P_{i} \setminus \{p '\}$. This permits us to linearize inequality (\ref{new_constraint_gamma}) with simple algebraic manipulations: we take the $\Gamma^\text{th}$-root of both sides of the inequality and then rely on the aforementioned property, giving us the following inequalities:
\begin{align}
    \sum_{p \in P_{i}}{\left(h_i^{1/\Gamma} - (h_i - p)^{1/\Gamma}\right) y_{ip}} + \sum_{j: (j,i) \in A}{d_{ji} z_{ji}} \geq h_{i}^{1/\Gamma} x_{i} & & \forall i \in V. \label{new_propagation_constraint_gamma_base} &
\end{align}

We also observe that there are cases where the incentive $p \in P_{i}$ offered to a node $i \in V$ is greater than the activation threshold for that node, that is, $p > h_{i}$. Since $h_{i} - p$ can be a negative value, it is desirable that $h_{i} - p$ be replaced by $0$ in these cases, leading to the following strengthened inequality:
\begin{align}
    \sum_{p \in P_{i}}{\left(h_i^{1/\Gamma} - \max (0, h_i - p)^{1/\Gamma}\right) y_{ip}} + \sum_{j: (j,i) \in A}{d_{ji} z_{ji}} \geq h_{i}^{1/\Gamma} x_{i} & & \forall i \in V. \label{new_propagation_constraint_gamma_max} &
\end{align}

Finally, we prove in Appendix \ref{appendix:proof_equivalent_constr} that Constraints \eqref{model:constraint_propagation} can be lifted into the following inequality:
\begin{align}
    \sum_{p \in P_{i}}{\left(\ceil{h_i^{1/\Gamma}} - \ceil{\max (0, h_i - p)^{1/\Gamma}}\right) y_{ip}} + \sum_{j: (j,i) \in A}{d_{ji} z_{ji}} \geq \ceil{h_{i}^{1/\Gamma}} x_{i} & & \forall i \in V. \label{new_propagation_constraint_gamma} &
\end{align}

Overall, this transformation models the original nonlinear function through new values of the incentives and associated costs in a linear context. Therefore, within the context of this formulation, any case with $\Gamma \neq 1$ can be reduced to the linear case with $\Gamma = 1$ using altered incentives and costs. Hence, for the sake of simplicity in the presentation of the formulations and algorithms, we will keep the original Constraint (\ref{model:constraint_propagation}) in the model from this point on, knowing that cases with $\Gamma \neq 1$ are handled by a transformation of the model parameters. 

\subsection{Influence Cover Cut (ICC)}
\label{subsec:cover_cut}

Having extended the breadth of applicability of the (ARC) formulation, we now focus on this formulation by proposing a new family of cutting planes for it, proving that these cuts are valid, and providing an efficient separation algorithm. This new cut will be called \emph{Influence Cover Cut (ICC)}, because it uses ideas present in cover cuts for the Knapsack Problem \citep{crowder1983solving}. This cut also shares many similarities with the generalized propagation constraints proposed in \citet{fischetti2018least} and discussed in Section~\ref{sec:Problem}, which are expressed on an extended formulation in which each variable corresponds to a set of influences that activates a given individual. The number of variables in this formulation grows exponentially with the number of neighbors of each individual.

In contrast, the proposed ICCs have the advantage of being expressed with a polynomial number of variables, namely the variables $y_{ip}$ describing the incentive $p$ offered to a node $i$, as well as the variables $x_{i}$ and $z_{ij}$ that describe which node is active and whether node $i$ influences a node $j$. This permits us to use this type of cover cuts within an (ARC) formulation.

The ICCs rely on the following logic: given a set of nodes $X \subseteq V$ containing an active node~$k$, at least one node $i \in X$ must receive an incentive $p$ and influences from nodes that are \textbf{not in} set $X$ up to a level that activates $i$. As such, node $i$ can initiate propagation within the set $X$, to activate the node $k$. In other words, there exists $U \subseteq V$ such that $f_{i}(U \setminus X, p) \geq h_{i}$. This leads to the following cuts.

\begin{proposition}
Consider a node $k$ and a subset of nodes $X \subseteq V$ such that $k \in X$. For each $i \in X$, select a subset of neighbors $\tilde{N}_{i} \subseteq N_{i} \setminus X$ and an incentive level $\tilde{p}_{i} \in P_{i}$ that are insufficient to activate the node $i$, i.e., such that:
\begin{eqnarray}
    \left(\sum_{j \in \tilde{N}_{i}}{d_{ji}}\right)^{\Gamma} + \tilde{p}_{i} < h_{i} \label{cover_cut:inequality1}
\end{eqnarray}
\noindent Then, the following inequality is valid:
\begin{eqnarray}
    \sum_{i \in X}{\sum_{\substack{
   p \in P_{i} \\
   p > \tilde{p}_{i}
}}{y_{ip}}} + \sum_{i \in X}{\sum_{j \in N_{i} \setminus (X \cup \tilde{N}_{i})}{z_{ji} \geq x_{k}.}} \label{cover_cut:inequality2}
\end{eqnarray}
\end{proposition}

\noindent
\proof{Proof.}
Let $(\bar{x}, \bar{y}, \bar{z})$ be a feasible solution. Given $k, X, \tilde{N}, \tilde{P}$ satisfying the conditions above and if $\bar{x}_{k} = 0$, then the constraint (\ref{cover_cut:inequality2}) is satisfied. Otherwise, let $k'$ be the first individual activated on $X$ (possibly $k' = k$). If $k'$ receives an incentive greater than $\tilde{p}_{k'} $, then there is a $p > \tilde{p}_{k'}$ such that $\bar{y}_{k'p} = 1$, which is sufficient to satisfy Constraint~(\ref{cover_cut:inequality2}).
    
In the remaining case, we demonstrate that (\ref{cover_cut:inequality2}) is satisfied by contradiction. Assume that this constraint is not satisfied and, therefore, that all influences received by $k'$ come from $X$ or $\tilde{N}_{k'}$. Since $k'$ is the first node activated in $X$, influences that come from $X$ cannot be used. Let $U = \left \{j \in N_{k'} \setminus X | \bar{z}_{jk'} = 1 \right \}$, and $p \in P_{k'} $ is the only value such that $\bar{y}_{k' p} = 1$. Since \eqref{cover_cut:inequality2} is not satisfied, $U \subseteq \tilde{N}_{k'}$ and $p \leq \tilde{p}_{k'}$. Thus, by inequality~\eqref{cover_cut:inequality1}, $f_{k'}(U, p) \leq \left (\sum_{j \in \tilde{N}_{k'}} {d_{j k'}} \right)^{\Gamma} \hspace*{-0.05cm} + \tilde{p}_{k'} < h_{k'}$, contradicting the assumption that $k'$ is activated.\Halmos
\endproof

\vspace*{1em}
\begin{example}
Consider an instance with $n=20$ nodes, average node degree $K=8$, minimum fraction of nodes to be activated $\alpha = 1.0$, $\Gamma = 1.0$ and set of possible incentives $P~=~ \{0.0, 9.0, 17.0, 26.0, 34.0\}$. Figure~\ref{fig:example_num_cover_cut} illustrates the cut at a given iteration of the branch-and-cut, in which the set $X = \{13, 19 \}$, $k = 13$, $\tilde{p}_{i} = 0$, $x_{13} = 1$ and $x_{19} = 1$. In this instance, the activation threshold for nodes $h_{13} = 9$, $h_{19} = 4$, as well as the list of incentives for nodes $13$ and $19$ is given by: $P_{13} = \{0.0, 9.0\}$ and $P_{19} = \{0.0, 9.0 \}$ respectively. Each arc $(i, j)$ of Figure~\ref{fig:example_num_cover_cut} represents the value $d_{ij}$, influence that the node $i$ has on the node $j$.

\begin{figure}[htbp]
    \centering
    \includegraphics[scale=0.7]{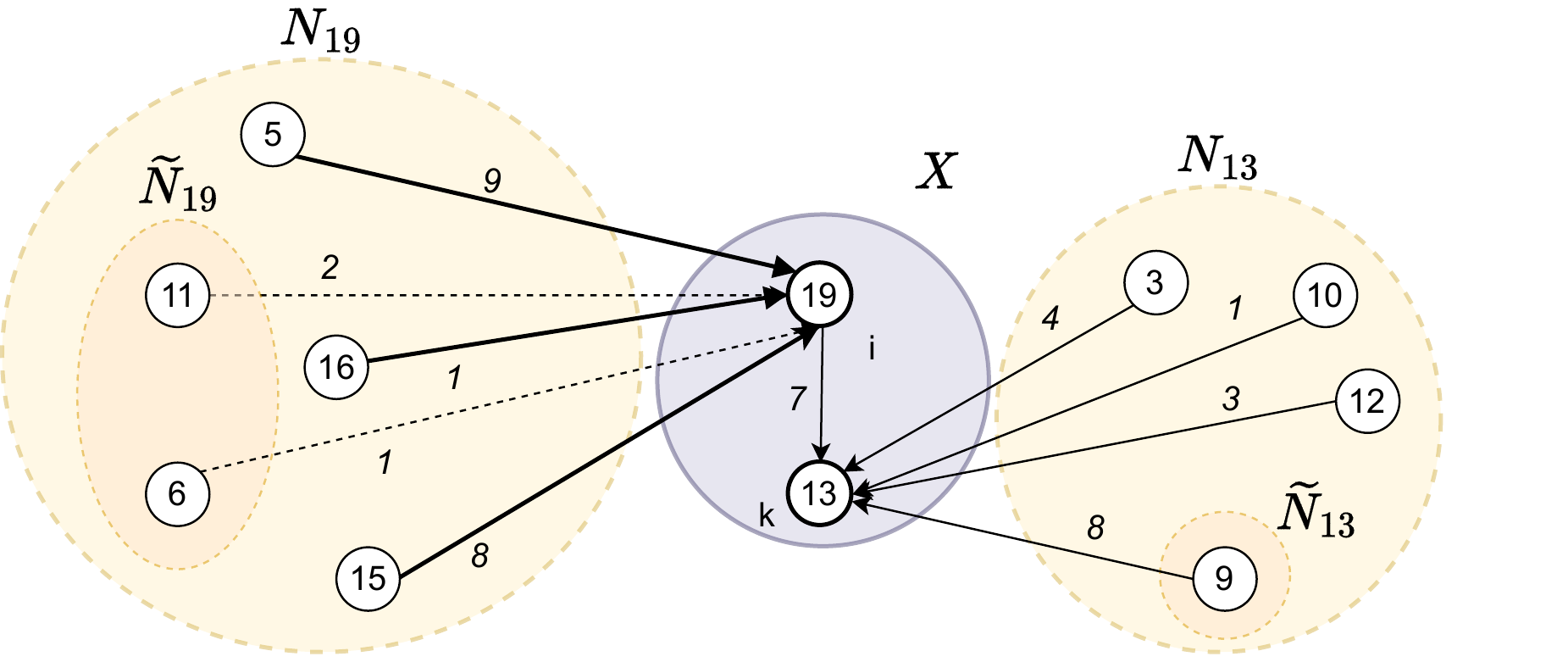}
    \caption{Illustration of the influence cover cut on a small instance}
    \label{fig:example_num_cover_cut}
\end{figure}

In this example, the set $\tilde{N}_{19}$ together with the incentive $\tilde{p}_{19}$ are not enough to activate the node $i = 19$. Thus, the following cut can be added to the model:
\begin{align}
    y_{13, 9} + y_{13, 17} + y_{19, 9} + y_{19, 17} + z_{3, 13} + &  \nonumber \\
    z_{10, 13} + z_{12, 13} + z_{5, 19} + z_{15, 19} + z_{16, 19} & \geq x_{13}.
     \label{cover_cut:inequality_example}
\end{align}

\end{example}

The incentive $\tilde{p}_{i}$ and the set of individuals outside the set $X$ that exert influences $\tilde{N}_{i}$ are chosen a priori as input parameters of the cut. They must be selected in order to satisfy~(\ref{cover_cut:inequality1}). Typically, in any solution of the linear relaxation, there will be enough incentives and influences outside set $X$ to activate the node $k$. Therefore, to locate a violated cut, we must remove some nodes from the set $\tilde{N}_{i}$ and reduce $\tilde{p}_{i}$ until the sum is not sufficient for activation.

We will refer to as ICC the model formed by Constraints \eqref{model:obj}--\eqref{model:constraint_z} with inequality \eqref{cover_cut:inequality2}. We now describe the separation algorithm designed to identify valid cuts.

\noindent
\subsubsection{Separation of the influence cover cuts.}
\label{subsection:icc_separation} \hfill

The separation problem of the inequalities \eqref{cover_cut:inequality2} is solved using an exact algorithm based on a mixed-integer programming (MIP) formulation, and the resulting cuts are added only at the root of the branch-and-cut tree. Let $(\bar{x}, \bar{y}, \bar{z})$ be a solution of the linear relaxation of model ARC. Our separation algorithm solves a MIP, described in the following, for each node $k \in V$ such that $\bar{x}_{k} > \epsilon$ (a predefined tolerance on the cut violation). Violations happen when the left side of Inequality (\ref{cover_cut:inequality2}) is smaller than the right side, so the separation MIP consists of a minimization problem on the left side to maximize the violation. A valid cut is obtained when the objective function of the MIP defined for node $k$ is smaller than $\bar{x}_{k} - \epsilon$. Thus, the algorithm can add up to one cut for each $k$ at each round.

The formulation of this MIP, which we will refer to as (SEP), uses five sets of variables. 
The variables $s_{i} ~ \forall i \in V$ indicate whether node $i$ belongs to the subset $X \subseteq V$ or not. The variables $y_{ip}^{0}$ and $z_{ji}^{0}$ represent respectively the incentive $\tilde{p}_{i}$ and nodes $j \in \tilde{N}_{i}$ that \textbf{are not} able to activate the node $i$. Moreover, the variables $y_{ip}^{1}$ and $z_{ji}^{1}$ represent the incentives $p > \tilde{p}_{i}$ and nodes $i \in X, j \in N_{i} \setminus(X \cup \tilde{N}_{i})$, and consequently, correspond to the variables $y_{ip}$ and $z_{ji}$ that are considered in the cut. In this way, the algorithm determines which nodes should be included in the set $\tilde{N}_{i}$ and incentive $\tilde{p}_{i}$ in order to maximize the violation of (\ref{cover_cut:inequality1}). The model is described below:

\noindent
\textbf{Variables:}
\begin{itemize}
    \item { $s_{i} = 1$ if $i \in X$ and $0$ otherwise}
    \item { $y_{ip}^{0} = 1$ if $p = \tilde{p}_{i}$ and $0$ otherwise}
    \item { $y_{ip}^{1} = 1$ if $i \in X$, $p > \tilde{p}_{i}$ and $0$ otherwise }
    \item { $z_{ji}^{0} = 1$ if $j \in \tilde{N}_{i}$ and $0$ otherwise}
    \item { $z_{ji}^{1} = 1$ if $i \in X$, $j \in  N_{i} \setminus (X \cup \tilde{N}_{i})$ and $0$ otherwise }
\end{itemize}

\noindent
\textbf{Formulation:}
\begin{align}[left = \hspace*{-1cm}\text{(SEP)}]
    \min \hspace*{0.4cm} & \sum_{i \in V}{\sum_{p \in P_{i}}{\bar{y}_{ip}y^{1}_{ip}}} + \sum_{i \in V}{\sum_{j \in N_{i}}{\bar{z}_{ji}z^{1}_{ji}}} \label{model:cover_cut_obj}\\
    \text{s.t.} \hspace*{0.4cm} & \sum_{p \in P_{i}}{p y_{ip}^{0}} + \sum_{j \in N_{i}}{d_{ji} z_{ji}^{0}} \leq h_{i} - 1 & \forall i \in V \label{model:cover_cut_cover} \\
    & y_{ip}^{1} \geq s_{i} - \sum_{\substack{q \in P_{i} \\ q \geq p}}{y_{iq}^{0}} & \forall i \in V, \forall p \in P_{i} \label{model:cover_cut_forcing_y1} \\
    & z_{ji}^{1} \geq s_{i} - s_{j} - z_{ji}^{0} & \forall i \in V, \forall j \in N_{i} \label{model:cover_cut_forcing_z1}\\
    & s_{k} = 1 & \label{model:cover_cut_k_in_X}\\
    & s_{i} \in \{0, 1\} & \forall i \in V \label{model:cover_cut_s}\\
    & z_{ji}^{1} \in \{0, 1\} & \forall i \in V, \forall j \in N_{i} \label{model:cover_cut_z1}\\
    & z_{ji}^{0} \in \{0, 1\} & \forall i \in V, \forall j \in N_{i} \label{model:cover_cut_z0}\\
    & y_{ip}^{1} \in \{0, 1\} & \forall i \in V, \forall p \in P_{i} \label{model:cover_cut_y1}\\
    & y_{ip}^{0} \in \{0, 1\} & \forall i \in V, \forall p \in P_{i}. \label{model:cover_cut_y0}
\end{align}

Objective~(\ref{model:cover_cut_obj}) maximizes the violation of the inequality (\ref{cover_cut:inequality2}). Constraints~(\ref{model:cover_cut_cover}) ensure that the set $\tilde{N}_{i}$ and an incentive $\tilde{p}_{i}$ are \textbf{not} able to activate a node $i \in V$, whereas inequalities~(\ref{model:cover_cut_forcing_y1})~and~(\ref{model:cover_cut_forcing_z1}) are used to link $y^{1}$ with $y^{0}$ and $z^{1}$ with $z^{0}$, respectively. Finally, Constraint~\eqref{model:cover_cut_k_in_X} ensures that node $k$ is in set $X$.

As seen in Section~\ref{sec:extension_gamma}, the case where $\Gamma \neq 1.0$ can be reduced to a case where $\Gamma = 1.0$ with different incentives and costs, and thus this separation routine is also applicable to the general version of the GLCIP. 

It is worth noting that both $\bar{y}$ and $\bar{z}$ have many zero components. Thus, there may be several optimal solutions to the separation MIP, some of them not using the maximum values of $\tilde{p}_i$ that satisfy \eqref{model:cover_cut_cover}. Similarly, the sets $\tilde{N}_i$ may not be maximal with respect to \eqref{model:cover_cut_cover}. We observe that directly using the cuts associated to these solutions substantially delays the convergence of the cut generation. Thus, we post-process the generated cuts by increasing the value of $\tilde{p}_i$ as much as possible, and then greedily complete $\tilde{N}_i$ by selecting the individuals with greater influence over $i$.

A \emph{cut round} is a run of the separation algorithm for all $k$, and as shown before, a cut round can add several cuts. The number of cutting rounds will be set as a parameter of our algorithm, along with a maximum time limit for cut separation. Therefore, the separation algorithm stops adding cuts when it achieves this time limit or when the algorithm reaches the defined number of cutting rounds, whichever condition is valid first. These two limits guarantee that the method does not spend excessive time on an individual run of the separation algorithm.

\subsection{A Special Case of the Influence Cover Cut}
\label{sec:icc+}

Let $X \subseteq V$, such that
\begin{eqnarray}
    |X| > (1 - \alpha) |V|. \label{eq:bigsubset}
\end{eqnarray}

Since there are not enough individuals outside the set $X$ to meet the minimum activation constraint, at least one individual in $X$ must be activated. Hence, for each $i \in X$, if the influences from $\tilde{N}_{i} \subseteq N_{i} \setminus X (i \in X)$ added to the incentive $\tilde{p}_{i} \in P_{i} (i \in X)$ are not enough to activate node $i$, then the following inequalities are also valid: 
\begin{eqnarray}
    \sum_{i \in X}{\sum_{\substack{
   p \in P_{i} \\
   p > \tilde{p}_{i}
}}{y_{ip}}} + \sum_{i \in X}{\sum_{j \in N_{i} \setminus (X \cup \tilde{N}_{i})}{z_{ji}}}\geq 1. \label{cover_cut:inequality3}
\end{eqnarray}

So, we can replace the right hand side $x_{i}$ with $1$ in inequalities \eqref{cover_cut:inequality2}. Thus, to guarantee that \eqref{eq:bigsubset} is satisfied in the separation MIP, it is necessary to replace the constraint \eqref{model:cover_cut_k_in_X} by:
\begin{eqnarray}
    \sum_{i \in V}{s_{i}} \geq \floor{(1 - \alpha) ~ |V|} + 1
\end{eqnarray}

Detecting and adding these new cuts is much simpler, since only one separation needs to be performed instead of $|V|$ in the case of the \emph{influence cover cuts}. This change significantly speeds up cut separation, making it possible to add these cuts in the entire branch-and-cut tree instead of limiting them to the root node.

We will refer to as ICC+ the model formed by Constraints \eqref{model:obj}--\eqref{model:constraint_z}, with inequalities \eqref{cover_cut:inequality2} and \eqref{cover_cut:inequality3}. 

\subsection{Lifted Influence Cover Cut}
\label{sec:lifted_icc+}

The cut described in Section~\ref{subsec:cover_cut} relies on the variables $z_{ji}$ to express the possible influence of external nodes of $i$. To further strengthen the model and make the separation algorithm lighter, we can express an inequality similar to (\ref{cover_cut:inequality2}) only in terms of the variables $x_{j}$ and $y_{ip}$ as follows. Consider a node $k$ and two subsets of nodes $X \subseteq V$ and $\tilde{N}_{X} \subseteq V \setminus X$ such that $k \in X$. Then, for each $i \in X$, select an incentive level $\tilde{p}_{i} \in P_{i}$ that, together with all incentives coming from set $\tilde{N}_{X}$, are insufficient to activate the node $i$. Namely, we have $\left(\sum_{j \in \tilde{N}_{X}} d_{ji}\right)^{\Gamma} + \tilde{p}_{i} < h_{i} \label{lift_cover_cut:inequality1}$ for all $i \in X$. Then, the following inequality is valid:

\begin{equation}
    \sum_{i \in X}{\sum_{\substack{
   p \in P_{i} \\
   p > \tilde{p}_{i}
}}{y_{ip}}} + \sum_{j \in V \setminus (X \cup \tilde{N}_{X})}{x_{j}} \geq x_{k}. \label{ineq:lifted_cover_cut}
\end{equation}

Consider, e.g., a set $X=\{1, 2\}$ and an outer node $3$ that influences the nodes $1$ and $2$. The original ICC cut would contain both variables $z_{13}$ and $z_{23}$ on the left side. Consequently, if node $3$ becomes activated it contributes twice to the sum, instead of once only when using inequality~\eqref{ineq:lifted_cover_cut}. Doing so also makes the separation lighter, since two-indexed variables are no longer necessary.

We empirically tested the combination of this cut with the ICC described in Section \ref{subsec:cover_cut} but observed that the model was computationally heavier. As a consequence, using only one of the two cuts seemed preferable.

\subsubsection{Separation of the lifted influence cover cuts.} \label{subsection:lifeted_cover_cut_sep}\hfill

The separation problem for inequalities \eqref{ineq:lifted_cover_cut} is also solved using an exact algorithm based on a MIP formulation, in a similar way as in Section \ref{subsection:icc_separation}. In preliminary analyses, this new cut yields good results, making the cut and separation algorithm much faster. This allowed us to add cuts of this type through the entire branch-and-cut tree, instead of just at the root as in the case of the ICC.

The MIP formulation of this separation problem uses five sets of variables. 
The variables~$s_{i}$ for $i \in V$ indicate whether node $i$ belongs to the subset $X \subseteq V$ or not. The variables $y_{ip}^{0}$ and $x_{j}^{0}$ represent respectively the incentive $\tilde{p}_{i}$ and nodes $j \in \tilde{N}_{X}$ that \textbf{are not} able to activate any node $i \in X$. Moreover, the variables $y_{ip}^{1}$ and $x_{j}^{1}$ represent the incentives $p > \tilde{p}_{i}$ and nodes $j \in V \setminus(X \cup \tilde{N}_{X})$, and consequently, correspond to the variables $y_{ip}$ and $x_{j}$ that are considered in the cut. In this way, the algorithm determines which nodes should be included in the set $\tilde{N}_{X}$ and incentive $\tilde{p}_{i}$ to maximize the violation of (\ref{ineq:lifted_cover_cut}).
The model is described below:

\noindent
\textbf{Variables:}
\begin{itemize}
    \item { $s_{i} = 1$ if $i \in X$ and $0$ otherwise}
    \item { $y_{ip}^{0} = 1$ if $p = \tilde{p}_{i}$ and $0$ otherwise}
    \item { $y_{ip}^{1} = 1$ if $i \in X$, $p > \tilde{p}_{i}$ and $0$ otherwise }
    \item { $x_{j}^{0} = 1$ if $j \in \tilde{N}_{X}$ and $0$ otherwise }
    \item { $x_{j}^{1} = 1$ if $j \in V \setminus(X \cup \tilde{N}_{X})$ and $0$ otherwise }
\end{itemize}

\noindent
\textbf{Formulation:}
\begin{align}[left = \hspace*{-1cm}\text{($SEP_{2}$)}]
    \min \hspace*{0.4cm} & \sum_{i \in V}{\sum_{p \in P_{i}}{\bar{y}_{ip}y^{1}_{ip}}} + \sum_{j \in V}{{\bar{x}_{j}x^{1}_{j}}} \label{model:lift_cover_cut_obj}\\
    \text{s.t.} \hspace*{0.4cm} & \sum_{p \in P_{i}}{p y_{ip}^{0}} + \sum_{j \in N_{i}}{d_{ji} x_{j}^{0}} + s_{i} \sum_{j \in N_{i}}{d_{ji}} \leq \sum_{j \in N_{i}}{d_{ji}} + h_{i} - 1 & \forall i \in V \label{model:lift_cover_cut_cover} \\
    & x_{i}^{1} \geq 1 - s_{i} - x_{i}^{0} & \forall i \in V \label{model:lift_cover_cut_forcing_x1} \\
    & y_{ip}^{1} \geq s_{i} - \sum_{\substack{q \in P_{i} \\ q \geq p}}{y_{iq}^{0}} & \forall i \in V, \forall p \in P_{i} \label{model:lift_cover_cut_forcing_y1} \\
    & s_{k} = 1 & \label{model:lift_cover_cut_k_in_X}\\
    & s_{i} \in \{0, 1\} & \forall i \in V \label{model:lift_cover_cut_s}\\
    & x_{i}^{1} \in \{0, 1\} &\forall i \in V \label{model:lift_cover_cut_x1}\\
    & x_{i}^{0} \in \{0, 1\} & \forall i \in V \label{model:lift_cover_cut_x0}\\
    & y_{ip}^{1} \in \{0, 1\} & \forall i \in V, \forall p \in P_{i} \label{model:lift_cover_cut_y1}\\
    & y_{ip}^{0} \in \{0, 1\} & \forall i \in V, \forall p \in P_{i}. \label{model:lift_cover_cut_y0}
\end{align}

Objective~\eqref{model:lift_cover_cut_obj} maximizes the violation of the inequality \eqref{ineq:lifted_cover_cut}. Constraints~\eqref{model:lift_cover_cut_cover} ensure that the set $\tilde{N}_{X}$ and an incentive $\tilde{p}_{i}$ are \textbf{not} able to activate a node $i \in V$, whereas inequalities~\eqref{model:lift_cover_cut_forcing_x1}~and~\eqref{model:lift_cover_cut_forcing_y1} are used to link $x^{1}$ with $x^{0}$ and $y^{1}$ with $y^{0}$, respectively.

\subsubsection{A special case of the lifted influence cover cut} \label{subsection:lifted_cover_cut_rhs_1}\hfill

In a similar way as described in Subsection \ref{sec:icc+} we can also generate an additional cut based on the inequality \eqref{ineq:lifted_cover_cut} given by: 
\begin{equation}
    \sum_{i \in X}{\sum_{\substack{
   p \in P_{i} \\
   p > \tilde{p}_{i}
}}{y_{ip}}} + \sum_{j \in V \setminus (X \cup \tilde{N}_{X})}{x_{j}} \geq 1, \label{ineq:lifted_cover_cut_rhs_1}
\end{equation}

\noindent
when \eqref{eq:bigsubset} is satisfied.
As in the other formulation, in the separation (MIP) algorithm of the inequalities \eqref{ineq:lifted_cover_cut_rhs_1} it is also necessary to replace the constraint \eqref{model:lift_cover_cut_k_in_X} by:
\begin{eqnarray}
    \sum_{i \in V}{s_{i}} \geq \floor{(1 - \alpha) ~ |V|} + 1
\end{eqnarray}

Thus, we refer to as LICC+ the model formed by the constraints \eqref{model:obj}--\eqref{model:constraint_z} with the inequalities \eqref{ineq:lifted_cover_cut} and \eqref{ineq:lifted_cover_cut_rhs_1}.

\subsection{Compact Formulation}\label{subsec:compact_formulation}

Inspired by the idea presented in Section~\ref{sec:lifted_icc+}, we finally developed a compact formulation that uses only the incentive variables $y_{ip}$. Let $\tilde{p}_{i}(X)$ be the smallest incentive that activates the node $i \in X$ given the influences of all individuals (activated or not) of the set $V \setminus X$. So, the compact formulation (CF) is described below:
\begin{align}[left = \hspace*{-1cm}\text{(CF)}]
    \min \hspace*{0.4cm} & \sum_{i \in V}{\sum_{p \in P_{i}}{w_{ip} y_{ip}}} \label{cf:obj_function}\\
    \text{s.a.} \hspace*{0.4cm} & \sum_{p \in P_{i}}{y_{ip}} = 1, & \forall i \in V \label{cf:incentive_node}\\
    & \sum_{i \in X}{\sum_{\substack{p \in P_{i} \\
    p \geq \tilde{p}_{i}(X)}}{y_{ip}}} \geq 1, & \forall X \subseteq V ~ | ~ |X| > \floor{(1-\alpha)|V|}  \label{cf:cut_inequality}\\
    & y_{i,p} \in \{0, 1\}, & \forall i \in V; \forall p \in P_{i} \label{cf:y_var}
\end{align}

Objective~\eqref{cf:obj_function} minimizes the total cost of the offered incentives. Constraints~\eqref{cf:incentive_node} guarantee that all nodes will receive an incentive (which can be zero). Constraints~\eqref{cf:cut_inequality} have an effect similar to the cut \eqref{cover_cut:inequality3} and guarantee that for any set $X$ whose complementary set has less than $\alpha |V|$ individuals, there must be sufficient influences and incentives on the set $X$ to activate at least one node. Consequently, these restrictions guarantee that at least $\alpha |V|$ nodes will be active at the end, as will be demonstrated in Section~\ref{sec:cf_correctness}.

There is an exponential number of inequalities \eqref{cf:cut_inequality}, and a separation algorithm is needed to identify and add violated constraints to the formulation. The algorithm developed for this purpose shares many similarities with the separation algorithm for ICCs. Moreover, in the case where the current solution is integral, a simple iterative algorithm permits to do the separation. In the rest of this section, we describe the separation algorithms and prove the validity of the formulation.

\subsubsection{Separation procedure for the integral case.}
\label{subsection:cf_integral_separation}\hfill

When the current solution is integral, a simple algorithm permits separating inequalities~\eqref{cf:cut_inequality} by simulating the propagation of influence in the graph. As seen in Algorithm~\ref{model_y:heuristic}, this algorithm initially calculates which nodes are already active due to the incentives (Lines 2--4). Then, each active node, in turn, propagates its influence to its neighbors, possibly activating new nodes in the process (Lines 5--12). The algorithm finishes as soon as each activated node has propagated its influence without activating new nodes. If the size of the set of activated nodes at the end of the process is not smaller than $\ceil{\alpha|V|}$, then we have a feasible solution. Otherwise, the set of non-activated nodes represents a set $X$ that characterizes a violated inequality~\eqref{cf:cut_inequality}.

\begin{figure*}[htbp]
\centering
\scalebox{0.83}
{
\begin{minipage}{1.2\textwidth}
\begin{algorithm}[H]
\SingleSpacedXII
\linespread{1.25}\selectfont
\caption{GetNonActivated($\bar{y}$)}\label{model_y:heuristic}
$A \leftarrow \varnothing$, $A_{\textsc{new}} \leftarrow \varnothing$\;
\For{$i \in V$}
{
    $R_{i} \leftarrow h_i - \sum_{p \in P_i}{p \bar{y}_{ip}}$
    \tcp*{$R_{i}$ is the amount of influence needed to activate $i$}
    \textbf{if} $R_{i} \leq 0$ \textbf{then} $A_{\textsc{new}} \leftarrow A_{\textsc{new}} \, \cup \,  \{i\}$ \tcp*[f]{node $i$ is activated due to initial incentives}
}
\While{$A_{\textsc{new}} \neq \varnothing$\tcp*[f]{propagate all new activations}}{
    $A_{\textsc{iter}} \leftarrow A_{\textsc{new}}$ \tcp*{nodes activated in the previous iteration} 
    $A_{\textsc{new}} \leftarrow \varnothing$ \tcp*{newly activated nodes}
    \For{\emph{$i \in V$ such that $R_{i} > 0$}}
    {
        \For{\emph{$j \in N_{i}$ such that $j \in A_{\textsc{iter}}$}}
        {
            $R_{i} \leftarrow R_{i} - d_{ji}$\;
            \textbf{if} $R_{i} \leq 0$ \textbf{then} $A_{\textsc{new}} \leftarrow A_{\textsc{new}} \, \cup \,  \{i\}$ \tcp*[f]{node $i$ becomes activated}
         }
    }
    $A \leftarrow A \cup A_{\textsc{iter}}$ \tcp*{store activated nodes}
}
\textbf{return} $V \setminus A$ \tcp*{return all non-activated nodes}
\end{algorithm}
\end{minipage}
}
\end{figure*}

\subsubsection{Separation procedure for the fractional case.}
\label{subsection:cf_frac_separation}\hfill

In the case where $\bar{y}$ is fractional, we can solve the separation problem through a MIP in a similar fashion as in Section \ref{subsec:cover_cut}.

This MIP, which we will refer to as (SEP$_{CF}$), uses four sets of variables. Variables $s_{i}^{1}, ~ \forall i \in V$ and $s_{i}^{0}, ~ \forall i \in V$ indicate that the node $i$ belongs to $X$ or that it belongs to $V \setminus X$ respectively. Variables $y_{ip}^{0}$ and $y_{ip}^{1}$ may always assume value $0$ in an optimal solution when $i \not\in X$. If $i \in X$, variables $y_{ip}^{0}$ indicate that $p$ is the maximum incentive that, together with the influences that come from outside of the set $X$, is not sufficient to activate node $i$ ($p < \tilde{p}_{i}(X)$). Finally,  variables $y_{ip}^{1}$ indicate that $p$ and the influences that come from outside of the set $X$ are together sufficient to activate node $i$ ($p \geq \tilde{p}_{i}(X)$). The model is described below:

\noindent
\textbf{Variables:}
\begin{itemize}
    \item $s_{j}^{0} = 1$ if $j \notin X$ and $0$ otherwise;
    \item $s_{i}^{1} = 1$ if $i \in X$ and $0$ otherwise;
    \item $y_{ip}^{0} = 1$ if $p < \tilde{p}_{i}$ and $0$ otherwise;
    \item $y_{ip}^{1} = 1$ if $i \in X$, $p \geq \tilde{p}_{i}$ and $0$ otherwise.
\end{itemize}

\noindent
\textbf{Formulation:}
Let $\bar{y}$ represent the relaxed values of the variables $y$. The formulation follows:

\begin{align}[left = \hspace*{-1cm}\text{(SEP$_{CF}$)}]
\allowdisplaybreaks
    \min \hspace*{0.4cm} & \sum_{i \in V}{\sum_{p \in P_{i}}{\bar{y}_{ip}y^{1}_{ip}}} \label{model_y:cover_cut_obj}\\
    \text{s.t.} \hspace*{0.4cm} & \sum_{p \in P_{i}}{p y_{ip}^{0}} + \sum_{j \in N_{i}}{d_{ji} (s_{j}^{0} + s_{i}^{1} - 1)} \leq  h_{i} - 1 & \forall i \in V \label{model_y:cover_cut_cover} \\
    & y_{ip}^{1} \geq s_{i}^{1} - \sum_{\substack{q \in P_{i} \\ q \geq p}}{y_{iq}^{0}} & \forall i \in V, \forall p \in P_{i} \label{model_y:cover_cut_forcing_y1} \\
    & s_{j}^{0} \geq 1 - s_{j}^{1} & \forall j \in V, \label{model_y:cover_cut_forcing_x1} \\
    & \sum_{i \in V}{s_{i}^{1}} \geq \floor{(1 - \alpha) |V|} + 1 \label{model_y:cover_cut_separation_x}\\
    & s_{i}^{1} \in \{0, 1\} & \forall i \in V \label{model_y:cover_cut_s1}\\
    & s_{j}^{0} \in \{0, 1\} &\forall j \in V \label{model_y:cover_cut_s0}\\
    & y_{ip}^{1} \in \{0, 1\} & \forall i \in V, \forall p \in P_{i} \label{model_y:cover_cut_y1}\\
    & y_{ip}^{0} \in \{0, 1\} & \forall i \in V, \forall p \in P_{i}. \label{model_y:cover_cut_y0}
\end{align}

Objective~\eqref{model_y:cover_cut_obj} maximizes the violation of  inequality \eqref{cf:cut_inequality}.
Constraints \eqref{model_y:cover_cut_cover} guarantee that any incentive $p < \tilde{p}_{i}(X)$ together with the influences that come from outside of the set $X$ (the sum of $d_{ji}$ such that $s_j^0 = s_i^1 = 1$) are not able to activate the node $i$ for all $i \in X$.
Constraints \eqref{model_y:cover_cut_forcing_y1} and \eqref{model_y:cover_cut_forcing_x1} link the variables $y_{ip}^{1}$ to $y_{ip}^{0}$, and $s_{j}^{0}$ to $s_{j}^{1}$, respectively.
Note that, since only one variable need to be set to $1$ on the right-hand side of constraint~\eqref{model_y:cover_cut_forcing_y1} to satisfy it, there is no incentive to set more than one variable $y_{ip}^{0}$ equal to $1$ for the same node~$i$.
Note also that, when $s_{i}^{1} = 0$ ($i \not\in X$), the second sum of \eqref{model_y:cover_cut_cover} is always non-positive. Moreover,  \eqref{model_y:cover_cut_forcing_y1} allows both $y_{ip}^{0}$ and $y_{ip}^{1}$ to be equal to $0$ for all $p \in P_i$. As a consequence, the vertex $i$ is not counted in the objective function if $i \not\in X$.
Constraint \eqref{model_y:cover_cut_separation_x} ensures that the set $X$ has more than $\floor{(1 - \alpha) |V|}$ individuals. 

\subsubsection{Correctness of formulation (CF).}
\label{sec:cf_correctness}

\begin{proposition} \label{theorem:cf_solution_feasible}
Any integer solution $\bar{y}$ that satisfies \eqref{cf:obj_function}--\eqref{cf:y_var} is feasible.
\end{proposition}

\proof{Proof.}
By contradiction. Suppose that $\bar{y}$ is an infeasible integer solution that satisfies all Constraints~\eqref{cf:cut_inequality}. Then, applying Algorithm 1, as $\bar{y}$ is infeasible, less than $\ceil{\alpha |V|}$ individuals will be activated. Let $X = V \setminus A$ be the set returned by Algorithm 1.

Since all individuals outside of set $X$, by construction, have been activated, we have in $X$ more than $|V| - \ceil{\alpha |V|} = |V| +\floor{- \alpha |V|} = \floor{(1 - \alpha)|V|}$ individuals. 

Therefore, $|X| > \floor{(1-\alpha)|V|}$, and $X$ defines a valid constraint of \eqref{cf:cut_inequality}.
Moreover, the algorithm only stops when no individual in $X$ can be activated by the incentives of $\bar{y}$ and the influences that come from outside of $X$.
Hence, the left-hand side of the constraint of \eqref{cf:cut_inequality} defined by $X$ is $0$.
Thus, this constraint is violated, contradicting the hypothesis. \Halmos
\endproof

\begin{proposition}
Any feasible integer solution $\bar{y}$ satisfies \eqref{cf:obj_function}--\eqref{cf:y_var}.
\end{proposition}
\proof{Proof.}
It is a direct consequence of the description of formulation (CF). \Halmos
\endproof

\section{Computational Experiments}
\label{sec:Experiments}

We conducted extensive numerical experiments to evaluate the performance of the different models and the impact of the cuts. This section reports the results of these analyses. All experiments were run on a single thread of an Intel Xeon E5-2620v4, $2.1$ GHz. The models were implemented using Julia\footnote{\url{https://julialang.org}}, JuMP v0.18.5 \footnote{\url{https://jump.dev/JuMP.jl/v0.18/}} \citep{DunningHuchetteLubin2017}, and IBM ILOG CPLEX $12.10$ was used to solve the MILP models with default parameters.

All the material (instances, solutions, source code and scripts) needed to reproduce these experiments is openly available at the following address: \url{https://github.com/vidalt/Influence-Optimization}.

\subsection{Instances and Numerical Conventions}

The instance set for the GLCIP was introduced by \citet{fischetti2018least}. For each node-set size $n \in \{50, 75, 100\}$, average node degree $K \in \{4,8,12,16\}$, and rewiring probability $\beta \in \{0.1, 0.3\}$, the authors generated five different graphs using the Watts-Strogatz method \citep{watts1998collective}, leading to a total of $90$ different graphs. Each graph is then declined into nine instances of the GLCIP by considering different values of $\alpha \in \{0.1, 0.5, 1.0\}$ and $\Gamma \in \{0.9, 1.0, 1.1\}$. Overall, there are $90 \times 9 = 810$ test instances.

In each of these instances, the set of incentives is defined as $P_{i} = \{0, 0.25\hat{h}, 0.5\hat{h}, 0.75\hat{h}, \hat{h} \}$ for each node $i$, where $\hat{h} = \max_{i \in V}{\{h_{i}\}}$ is the highest activation threshold among all nodes, and the cost for offering incentive $p$ to node $i$ is calculated as $w_{ip} = p^{0.9}$.

As indicated to us by the authors, to achieve the same experimental conventions as in \citet{fischetti2018least}, it is necessary to round up the incentive values in $P_{i}$. Moreover, the incentive costs $w_{i,p}$ should be truncated, and the activation function $f_{i}(U, p)$ should be rounded to the nearest integer. The authors also identified a minor bug in their original implementation that impacted the results for $\Gamma = 0.9$. Therefore, they provided us with a corrected set of results for each instance to permit detailed comparisons. All the solutions are provided in the Github repository associated with this study.

In previous work, \citet{fischetti2018least} reported results for six different formulations named A, C, C$^{+}$, C$^{+}_{e}$, P$^{+}$, and P$^{ +}_{e}$. The first four models are used to design exact methods, whereas the last two are used for heuristic search. We will compare our models with A, C, and C$^{+}_{e}$ from \citet{fischetti2018least}. Model A serves as a baseline since our branch-and-cut method builds upon it by adding the proposed \emph{influence cover cuts}.
Included formulation C was selected, as it achieves the current best results for $\Gamma=0.9$. Finally, C$^{+}_{e}$ was selected since it achieves the best overall exact results and is particularly effective for $\Gamma \in \{1.0, 1.1\}$. We compare these algorithms with the following new formulations:
\begin{itemize}
    \item ICC: Model formed by Constraints \eqref{model:obj}--\eqref{model:constraint_z} with inequality \eqref{cover_cut:inequality2} as presented in Section~\ref{subsec:cover_cut}. We set a limit of $200$ cut rounds or a maximum separation time of $5$ minutes, and only separate cuts at the root node of the branch-and-cut tree.
    \item ICC+: Variation of model ICC which includes inequality \eqref{cover_cut:inequality3} as presented Section~\ref{sec:icc+}. The separation of inequality \eqref{cover_cut:inequality2} is the same as in model ICC. In contrast, the separation of inequality \eqref{cover_cut:inequality3} is pursued through the entire solution process.
    \item LICC+: Model formed by Constraints \eqref{model:obj}--\eqref{model:constraint_z} along with inequalities \eqref{ineq:lifted_cover_cut} and \eqref{ineq:lifted_cover_cut_rhs_1}, as presented in Section~\ref{subsection:lifted_cover_cut_rhs_1}. The separation algorithm for inequality \eqref{ineq:lifted_cover_cut} is used through the entire solution process, and the separation algorithm for inequality \eqref{ineq:lifted_cover_cut_rhs_1} is used as long as optimality gap $\geq 40\%$.
    \item CF: Compact formulation formed by Constraints \eqref{cf:obj_function}--\eqref{cf:y_var}, as presented in Section \ref{subsec:compact_formulation}. The separation algorithm for the fractional case of inequality \eqref{cf:cut_inequality} is used as long as optimality gap $\geq 10\%$.
\end{itemize}

In all the experiments, we use as cutoff the value of the best known exact solution reported in \citet{fischetti2018least} increased by one unit, and we set a time limit of 2 hours for all the algorithms. In the remainder of this section, we compare the different algorithms in terms of their number of instances solved in optimality, their final optimality gaps, and total execution time. Additional detailed results are provided in~Appendix \ref{appendix:detailed_results} and in the Github repository.

\subsection{Number of Optimal Solutions} 
\label{sec:number_optimal}

\begin{figure}[htpb!] 
    \hspace*{-1.7cm}
    \includegraphics[scale=0.6]{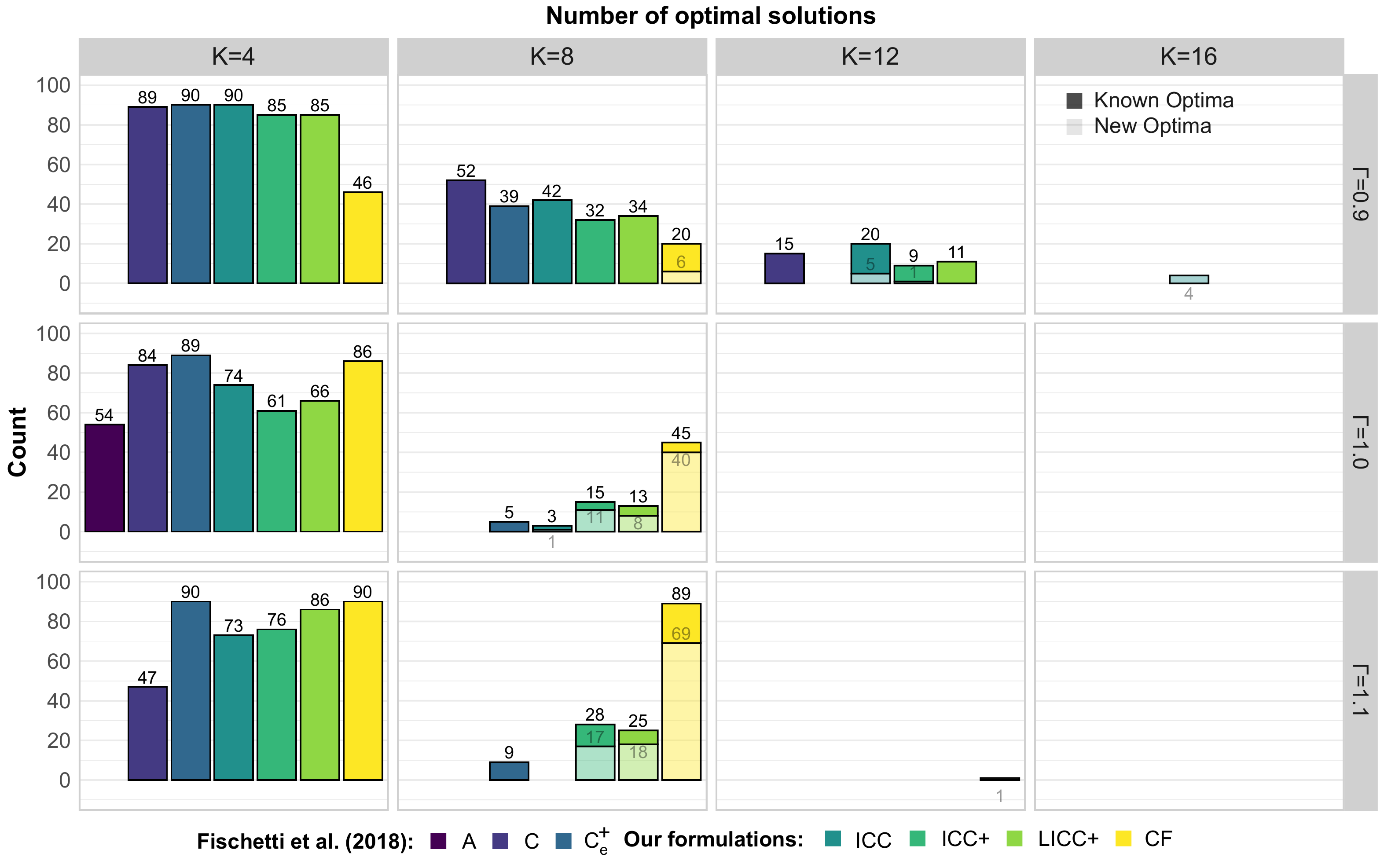}
    \caption{Nb. of instances solved to optimality by each formulation, for $\boldsymbol{K \in \{4, 8, 12, 16\}}$ and $\boldsymbol{\Gamma \in \{0.9, 1.0, 1.1\}}$}
    \label{fig:number_instances_solved}
\end{figure}

Figure~\ref{fig:number_instances_solved} compares the number of instances solved to optimality for each formulation. These results are presented as a faceted grid of bar plots, including a row for each value of $\Gamma \in \{0.9, 1.0, 1.1\}$ and a column for each value of $K \in \{4, 8, 12, 16\}$. The value (count) of each bar represents the number of instances optimally solved, with the darkest bars being the absolute number of optimal instances. In contrast, the light bars correspond to instances that were optimally solved but for which no optimal solution was known until now. Formulation~$A$ of \citet{fischetti2018least} is only applicable when $\Gamma = 1.0$, therefore the results of this formulation only appear in the middle group.

As seen in these results, the value of $\Gamma$ has a big impact on the number of solved instances, but this impact is different for each formulation. The proposed formulation CF, for example, achieves the smallest number of optimal solutions among all formulations when $\Gamma = 0.9$, while formulations C and ICC are achieve the largest number. In contrast, when $\Gamma \in \{1.0, 1.1\}$, the proposed CF formulation solves the largest number of instances, that is, four times the number of instances solved by the C formulation, and nearly twice the number of instances solved by the C$^{+}_{e}$ formulation of \citet{fischetti2018least}. Overall, it appears that the best formulation choice depends on the specific characteristics of the influence functions.

In contrast, the average degree in the network represented by parameter $K=8$ has a more consistent effect on the performance of all the formulations. As the graph gets denser, the number of possible combinations of influences able to activate a node dramatically increase, and solving the resulting GLCIP becomes increasingly difficult.

It is noteworthy that the proposed formulations could solve many instances for which the optimum was unknown so far: 125 new optimal solutions have been found. 
It is also noteworthy that formulation ICC could solve, for the first time, four instances with a degree $K=16$. Many of the newly solved instances have a larger degree (e.g., $K=8$) and involve supermodular influence functions with $\Gamma=1.1$, in which the marginal influence of each incentive grows with the number of incentives. For this combination of factors, previous formulations were fairly ineffective, but the CF formulation performs generally well, therefore this approach fills an important methodological gap.
 
\subsection{Optimality Gaps}
\label{sec:distribution_gaps}

\begin{figure}[htpb!]
    \centering
    \hspace*{-1.5cm}
    \includegraphics[scale=0.38]{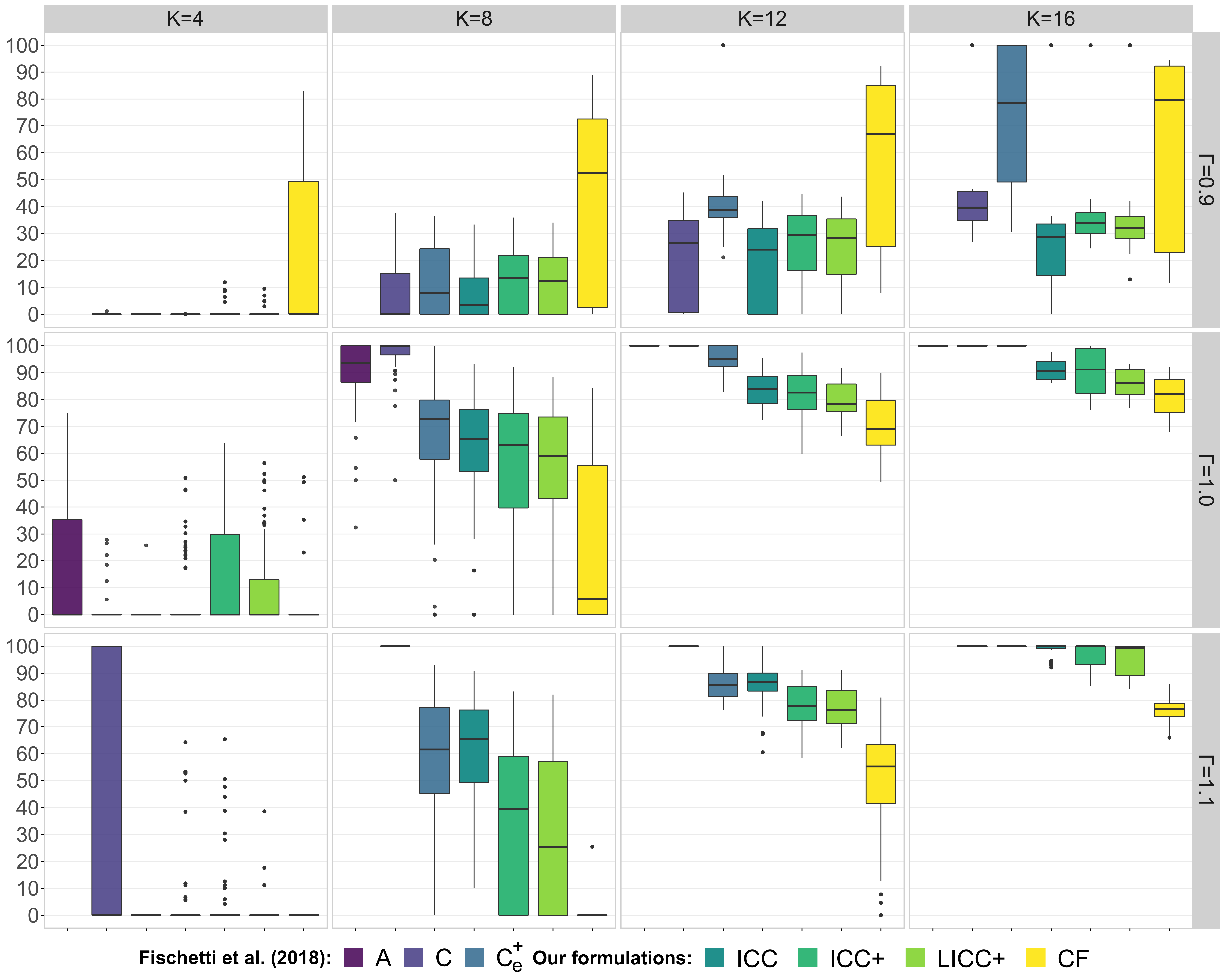}
    \caption{Final optimality gaps of the different formulations, for $\boldsymbol{K \in \{4, 8, 12, 16\}}$ and $\boldsymbol{\Gamma \in \{0.9, 1.0, 1.1\}}$}
    \label{fig:distribution_gap}
\end{figure}

If the branch-and-cut reaches the 2-hours time limit without identifying and proving the optimality of a solution, then it terminates with an optimality gap calculated as $\text{Gap}(\%) = 100 \times (Z_\textsc{ub} - Z_\textsc{lb}) / Z_\textsc{ub}$, where $Z_\textsc{ub}$ and $Z_\textsc{lb}$ are the values of the best upper and lower bounds found during the search. Figure \ref{fig:distribution_gap} represents the distribution of the optimality gaps achieved by the formulations at the end of the solution process. As previously, the results are presented on a faceted grid of boxplots. each row corresponds to a different value of $\Gamma~\in~\{0.9, 1.0, 1.1\}$ and each column corresponds to a different value of $K~\in~\{4, 8, 12.16\}$ (average node degree). Each boxplot within the graph is therefore associated with $90$ observations, and the whiskers extend to $1.5$ times the interquartile range.

As already noted in \citet{fischetti2018least}, we observe that the complexity of the instances (and therefore the final gap of the methods) increases with $K$. First, we can highlight that, just for $\Gamma=0.9$, the proposed ICC formulation obtains the smallest average gap over all formulations. In contrast, for $\Gamma \in \{1.0, 1.1\}$, our compact formulation CF achieves the best performance. It is also worth mentioning that only for $\Gamma=0.9$, the best results obtained by \citet{fischetti2018least} are achieved by the C formulation, while for $\Gamma \in \{1.0, 1.1\}$, they are achieved by the formulation~C$^{+}_{e}$.

In addition, we can also observe that smaller values of $\Gamma$ generally correlate with smaller gaps. This holds for all models except for CF, where this behavior appears to be reversed.

\subsection{Computational Effort}
\label{sec:distribution_execution_time}

\begin{figure}[htpb!]
    \hspace{-1.5cm}
    \includegraphics[scale=0.37]{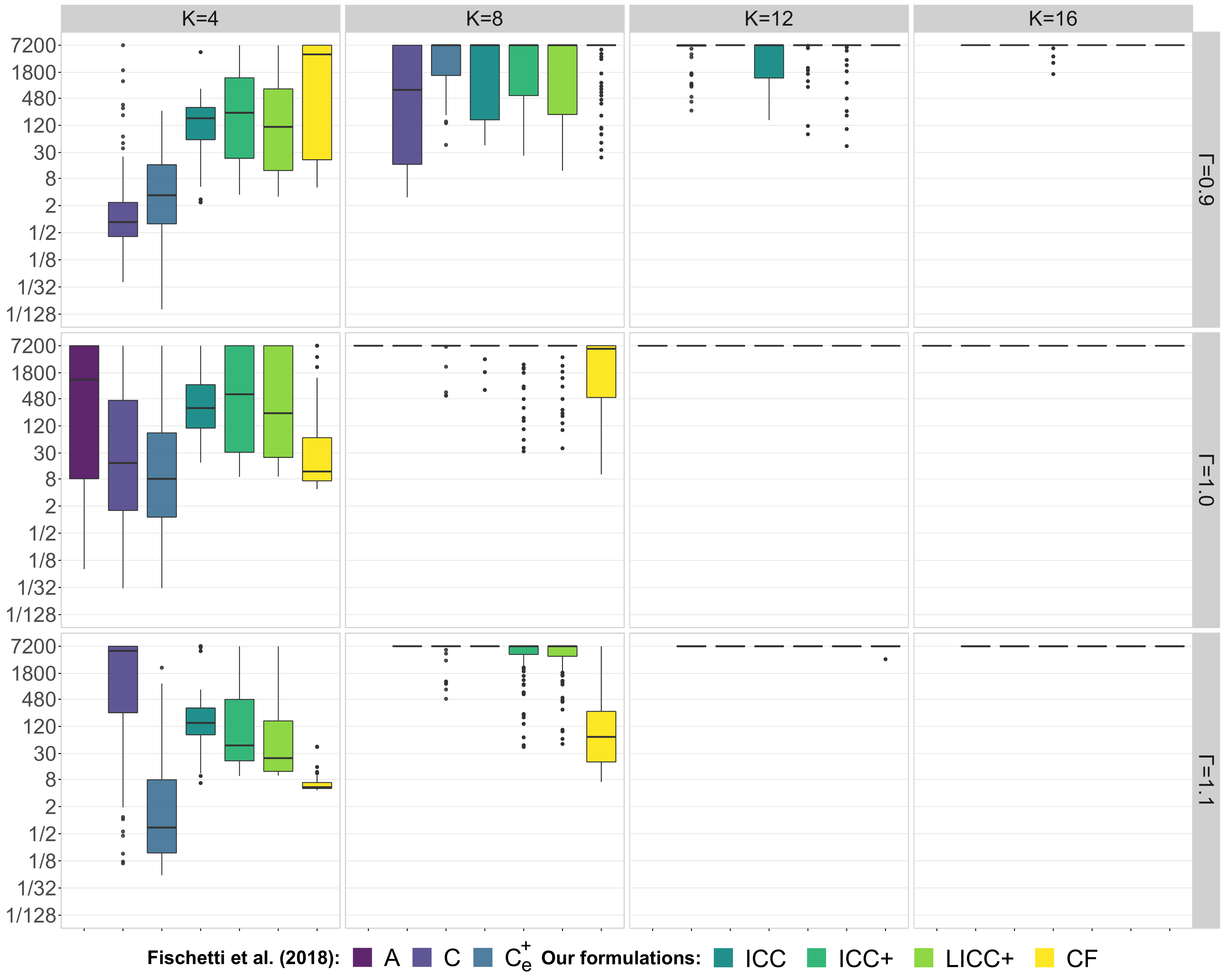}
    \caption{Computational time of the different formulations, for $\boldsymbol{K \in \{4, 8, 12, 16\}}$ and $\boldsymbol{\Gamma \in \{0.9, 1.0, 1.1\}}$}
    \label{fig:distribution_time}
\end{figure}

Finally, Figure~\ref{fig:distribution_time} compares the execution time of each formulation for the different values of $\Gamma$ and $K$. These results are presented using a logarithmic scale as a faceted grid of boxplots: one row for each value of $\Gamma~\in~\{0.9, 1.0, 1.1\}$, and one column for each value of $K~\in~\{4, 8, 12.16\}$.

As noted previously, the higher the value of $K$, the more difficult the instance becomes, and consequently, the execution time distribution tends to be closer to the 2-hour (i.e., 7200 seconds) limit. For $K = 4$, we observe that the C$^{+}_{e}$ formulation yields the shortest execution times, except for the case of $\Gamma=0.9$ where the C formulation is faster. For $K = 8$, the CF formulation achieves the shortest execution times, except for the case with $\Gamma=0.9$. Finally, for $K \in \{12, 16\}$ all formulations tend to attain the time limit, regardless of the value of~$\Gamma$.

\section{Conclusions}
\label{sec:Conclusions}

In this paper, we have contributed to better solving the GLCIP through new mathematical formulations and valid inequalities. We showed how to extend the arc-flow model presented in \citet{fischetti2018least} to address nonlinear cases ($\Gamma \neq 1$) not initially covered by it. Furthermore, we introduced new variations of this base formulation through new cutting planes and separation algorithms. Finally, we introduced a new compact formulation (CF) that is much simpler than the other existing ones and requires an adaptation of the valid inequalities and separation algorithms. As seen in our experimental results, this compact formulation permits solving over $120$ previously-open instances to proven optimality. In other cases, it also found the optimal values previously reported by \citet{fischetti2018least} or led to better bounds.

The research perspectives connected to this work are numerous.
They can concern the development of new pre-processing procedures, formulations, and stronger cuts, to further strengthen the model or make the separation algorithm lighter, with the effect of improving the algorithm performance to deal with larger (or denser) instances.
We also suggest further extending existing influence propagation models and methods to other variants of importance. For example, most studies assume that the influence probabilities are known in advance. Yet, data uncertainty is important in practice and should be taken into account. Moreover, whereas exact methods permit obtaining optimal solutions for networks of relatively small scale and act as a reliable benchmark, real applications of influence propagation often occur on huge networks. In this case, it is necessary to turn towards heuristic or metaheuristic approaches for influence optimization, and further developments along this line are likely to be needed.

Last but not least, influence optimization is one of several possible objectives. In practice, mitigation actions may be even more critically needed than propagation strategies, for example, to help preventing the spread of epidemics in populations or fake news on social networks. Developing effective mitigation measures is a critical and difficult research subject. One way to achieve this could be through a two-player game model in which one player seeks to maximize influence propagation, whereas the other player attempts to prevent or mitigate the spread. Given how hard it already is to formulate and solve the influence optimization problem in a general setting with a single decision maker, it is likely that such adversarial settings will remain a significant challenge in the years to come. We therefore hope that the work that we conducted on the simpler setting of influence optimization will permit progress in this direction.

\newpage
\ACKNOWLEDGMENT{%
This research has been partially funded by CAPES [Finance Code 001], CNPq [grant number 308528/2018-2], and FAPERJ [grant number E-26/202.790/2019] in Brazil. This support is gratefully acknowledged.
}



%
%
%
\begin{APPENDICES}
\newpage
\section{Proofs}
\label{appendix:proof_equivalent_constr}

\begin{proposition} \label{theorem:equivalent_constraints}
Constraints \eqref{new_constraint_gamma} and \eqref{new_propagation_constraint_gamma} are equivalent.
\end{proposition}

The proof will be broken down into two cases, depending on the values of $x_{i} $ and $y_{i, p}$.

\proof{Proof.}
First, if we assume that $x_{i} = 0$, then we have:
\begin{align}
\sum_{p \in P_{i}}{(\ceil{h_i^{1/\Gamma}} - \ceil{\max (0, h_i - p)^{1/\Gamma}}) y_{ip}} + \sum_{j: (j,i) \in A}{d_{ji} z_{ji}} \geq 0 & & \forall i \in V. \label{new_propagation_constraint_gamma2} &
\end{align}

Therefore, the lefthand side is always positive and both constraints \eqref{new_constraint_gamma} and \eqref{new_propagation_constraint_gamma} are satisfied.\\

Otherwise, we can assume that $x_{i} = 1$. A single incentive is offered to node $i$, therefore only one variable $y_{i, p}, \, \forall i \in V, p \in P_{i}$ is equal to~1.
Let $p'$ be the incentive offered at node $i$, i.e., we have $y_{i, p'} = 1$ and $y_{i, p} = 0, ~ \forall p \in P_{i} \setminus \{p'\}$. As a consequence, 
\begin{align}
\text{Constraint}~\eqref{new_propagation_constraint_gamma} & \Longleftrightarrow \ceil{h_{i}^{1/\Gamma}} - \ceil{\max (0, h_i - p')^{1/\Gamma}} + \sum_{j: (j,i) \in A}{d_{ji} z_{ji}} \geq \ceil{h_{i}^{1/\Gamma}} & & \forall i \in V \label{new_propagation_constraint_gamma3} & \\
& \Longleftrightarrow  \sum_{j: (j,i) \in A}{d_{ji} z_{ji}} \geq \ceil{\max (0, h_i - p')^{1/\Gamma}} & & \forall i \in V. \label{new_propagation_constraint_gamma4} &
\end{align}

We now need to separate two possible cases.
If $h_{i} - p'\leq 0$, then Constraint~\eqref{new_propagation_constraint_gamma4} is equivalent to $\sum_{j: (j,i) \in A}{d_{ji} z_{ji}} \geq 0$, which is satisfied regardless of the values of the variables $z$, as well as \eqref{new_constraint_gamma}.
Otherwise, we are in a situation where $h_{i} - p'> 0$, i.e., incentive $p'$ is insufficient to activate the node $i$. Given that $d_{ji}$ take integer values, we have:
\begin{align}
 \sum_{j: (j,i) \in A}{d_{ji} z_{ji}} \geq \ceil{\max (0, h_i - p')^{1/\Gamma}} &\Longleftrightarrow  \sum_{j: (j,i) \in A}{d_{ji} z_{ji}} \geq \ceil{ (h_i - p')^{1/\Gamma}}, & & \forall i \in V \label{new_propagation_constraint_gamma5} & \\
 & \Longleftrightarrow  \sum_{j: (j,i) \in A}{d_{ji} z_{ji}} \geq (h_i - p')^{1/\Gamma}, & & \forall i \in V \label{new_propagation_constraint_gamma6} &
 \end{align}
Raising both sides of this inequality to the power of $\Gamma$, we have:
\begin{align}
\sum_{j: (j,i) \in A}{d_{ji} z_{ji}} \geq (h_i - p')^{1/\Gamma} & \Longleftrightarrow  \left(\sum_{j: (j,i) \in A}{d_{ji} z_{ji}}\right)^{\Gamma} \geq ((h_i - p')^{1/\Gamma})^{\Gamma}, & & \forall i \in V \label{new_propagation_constraint_gamma7} & \\
& \Longleftrightarrow  \left(\sum_{j: (j,i) \in A}{d_{ji} z_{ji}}\right)^{\Gamma} \geq h_i - p', & & \forall i \in V \label{new_propagation_constraint_gamma8} & \\
& \Longleftrightarrow  p' + \left(\sum_{j: (j,i) \in A}{d_{ji} z_{ji}}\right)^{\Gamma} \geq h_i & & \forall i \in V. \label{new_propagation_constraint_gamma9} &
\end{align}
As only one incentive is paid for a node, we can replace $p'$ again with $\sum_{p \in P_{i}} {p y_{i, p}} $, obtaining the following inequality equivalent to \eqref{new_constraint_gamma}:
\begin{align}
\Longleftrightarrow  \sum_{p \in P_{i}}{p y_{i,p}} + \left(\sum_{j: (j,i) \in A}{d_{ji} z_{ji}}\right)^{\Gamma} \geq h_i, & & \forall i \in V. \label{new_propagation_constraint_gamma10} &
\end{align}
As a consequence, the equivalence of the inequalities has been proven in all cases. \Halmos
\endproof

\newpage
\section{Detailed results} 
\label{appendix:detailed_results}

This section presents finer-grained results which are aggregated (averaged) over the five instances for each combination of the $n, K, \beta$ and $\alpha$ factors. These results are provided in Tables~\ref{tab:results_avg_0_9}, \ref{tab:results_avg_1_0}, and~\ref{tab:results_avg_1_1} for $\Gamma \in \{0.9, 1.0, 1.1\}$ respectively. The best values are highlighted in bold.


\newpage
\begin{table}[H]
\centering
\scalebox{0.585}{
\begin{tabular}{rrrrrrrrrrrrrrrrrr}
\toprule
\multicolumn{1}{c}{\textbf{}} & \multicolumn{1}{c}{\textbf{}} & \multicolumn{1}{c}{\textbf{}} & \multicolumn{1}{c}{\textbf{}} & \multicolumn{1}{c}{} & \multicolumn{6}{c}{\textbf{Avg. gap (\%)}} & \multicolumn{1}{c}{} & \multicolumn{6}{c}{\textbf{Avg. total time (s)}} \\ \cmidrule(lr){6-11} \cmidrule(l){13-18} 
\multicolumn{1}{c}{\textbf{n}} & \multicolumn{1}{c}{\textbf{K}} & \multicolumn{1}{c}{$\bm{\beta}$} & \multicolumn{1}{c}{$\bm{\alpha}$} & \multicolumn{1}{c}{} & \multicolumn{1}{c}{$\bm{\mathrm{C}}$} & \multicolumn{1}{c}{$\bm{\mathrm{C}}^{+}_{e}$} & \multicolumn{1}{c}{$\bm{\mathrm{ICC}}$} & \multicolumn{1}{c}{$\bm{\mathrm{ICC+}}$} & \multicolumn{1}{c}{$\bm{\mathrm{LICC+}}$ } & \multicolumn{1}{c}{$\bm{\mathrm{CF}}$} & \multicolumn{1}{c}{} & \multicolumn{1}{c}{$\bm{\mathrm{C}}$} & \multicolumn{1}{c}{$\bm{\mathrm{C}}^{+}_{e}$} & \multicolumn{1}{c}{$\bm{\mathrm{ICC}}$} & \multicolumn{1}{c}{$\bm{\mathrm{ICC+}}$} & \multicolumn{1}{c}{$\bm{\mathrm{LICC+}}$} & \multicolumn{1}{c}{$\bm{\mathrm{CF}}$} \\ \midrule
\rowcolor[HTML]{EFEFEF} 
50 & 4 & 0.1 & 0.1 &  & \textbf{0.00} & \textbf{0.00} & \textbf{0.00} & \textbf{0.00} & \textbf{0.00} & \textbf{0.00} &  & \textbf{0.24} & 2.14 & 81.38 & 16.28 & 9.19 & 7.97 \\
\rowcolor[HTML]{EFEFEF} 
50 & 4 & 0.1 & 0.5 &  & \textbf{0.00} & \textbf{0.00} & \textbf{0.00} & \textbf{0.00} & \textbf{0.00} & 36.89 &  & \textbf{0.72} & 12.79 & 216.07 & 616.86 & 695.89 & 7200.00 \\
\rowcolor[HTML]{EFEFEF} 
50 & 4 & 0.1 & 1.0 &  & \textbf{0.00} & \textbf{0.00} & \textbf{0.00} & \textbf{0.00} & \textbf{0.00} & \textbf{0.00} &  & 6.19 & \textbf{1.00} & 88.04 & 153.10 & 57.53 & 16.18 \\
\rowcolor[HTML]{EFEFEF} 
50 & 4 & 0.3 & 0.1 &  & \textbf{0.00} & \textbf{0.00} & \textbf{0.00} & \textbf{0.00} & \textbf{0.00} & \textbf{0.00} &  & \textbf{0.16} & 3.59 & 30.53 & 9.01 & 9.81 & 7.46 \\
\rowcolor[HTML]{EFEFEF} 
50 & 4 & 0.3 & 0.5 &  & \textbf{0.00} & \textbf{0.00} & \textbf{0.00} & \textbf{0.00} & \textbf{0.00} & 40.88 &  & \textbf{0.64} & 13.63 & 144.94 & 681.03 & 1003.96 & 7200.00 \\
\rowcolor[HTML]{EFEFEF} 
50 & 4 & 0.3 & 1.0 &  & \textbf{0.00} & \textbf{0.00} & \textbf{0.00} & \textbf{0.00} & \textbf{0.00} & \textbf{0.00} &  & 15.82 & \textbf{5.99} & 241.47 & 719.27 & 208.80 & 26.83 \\
50 & 8 & 0.1 & 0.1 &  & \textbf{0.00} & \textbf{0.00} & \textbf{0.00} & \textbf{0.00} & \textbf{0.00} & \textbf{0.00} &  & \textbf{5.35} & 354.50 & 124.74 & 49.63 & 31.74 & 262.11 \\
50 & 8 & 0.1 & 0.5 &  & \textbf{0.00} & \textbf{0.00} & \textbf{0.00} & 11.20 & 7.77 & 64.94 &  & \textbf{18.93} & 2163.95 & 338.81 & 6037.68 & 5282.11 & 7200.00 \\
50 & 8 & 0.1 & 1.0 &  & 13.37 & 1.27 & 8.03 & 6.96 & 6.62 & \textbf{0.96} &  & 7200.00 & \textbf{2070.42} & 7200.00 & 6163.32 & 5707.16 & 3278.99 \\
50 & 8 & 0.3 & 0.1 &  & \textbf{0.00} & \textbf{0.00} & \textbf{0.00} & \textbf{0.00} & \textbf{0.00} & \textbf{0.00} &  & \textbf{4.41} & 342.35 & 63.78 & 28.28 & 22.94 & 140.96 \\
50 & 8 & 0.3 & 0.5 &  & \textbf{0.00} & 16.10 & 1.47 & 25.48 & 23.46 & 72.64 &  & \textbf{557.31} & 7200.00 & 2645.17 & 7200.00 & 7200.00 & 7200.00 \\
50 & 8 & 0.3 & 1.0 &  & 17.12 & 9.26 & 14.75 & 14.69 & 14.23 & \textbf{2.35} &  & 7200.00 & 7200.00 & 7200.00 & 7200.00 & 7200.00 & \textbf{4353.84} \\
\rowcolor[HTML]{EFEFEF} 
75 & 4 & 0.1 & 0.1 &  & \textbf{0.00} & \textbf{0.00} & \textbf{0.00} & \textbf{0.00} & \textbf{0.00} & \textbf{0.00} &  & \textbf{0.59} & 2.87 & 38.07 & 22.59 & 14.30 & 30.92 \\
\rowcolor[HTML]{EFEFEF} 
75 & 4 & 0.1 & 0.5 &  & \textbf{0.00} & \textbf{0.00} & \textbf{0.00} & \textbf{0.00} & \textbf{0.00} & 63.66 &  & \textbf{0.74} & 9.69 & 199.92 & 632.13 & 471.33 & 7200.00 \\
\rowcolor[HTML]{EFEFEF} 
75 & 4 & 0.1 & 1.0 &  & \textbf{0.00} & \textbf{0.00} & \textbf{0.00} & \textbf{0.00} & \textbf{0.00} & 0.91 &  & 73.23 & \textbf{3.42} & 203.60 & 1049.85 & 375.55 & 4056.09 \\
\rowcolor[HTML]{EFEFEF} 
75 & 4 & 0.3 & 0.1 &  & \textbf{0.00} & \textbf{0.00} & \textbf{0.00} & \textbf{0.00} & \textbf{0.00} & \textbf{0.00} &  & \textbf{0.35} & 10.66 & 39.08 & 17.88 & 10.04 & 122.34 \\
\rowcolor[HTML]{EFEFEF} 
75 & 4 & 0.3 & 0.5 &  & \textbf{0.00} & \textbf{0.00} & \textbf{0.00} & 5.07 & 3.84 & 71.43 &  & \textbf{2.97} & 109.58 & 226.63 & 5375.22 & 5224.02 & 7200.00 \\
\rowcolor[HTML]{EFEFEF} 
75 & 4 & 0.3 & 1.0 &  & \textbf{0.00} & \textbf{0.00} & \textbf{0.00} & \textbf{0.00} & \textbf{0.00} & 2.16 &  & 15.92 & \textbf{4.95} & 212.49 & 2541.12 & 1338.51 & 3535.69 \\
75 & 8 & 0.1 & 0.1 &  & \textbf{0.00} & \textbf{0.00} & \textbf{0.00} & \textbf{0.00} & \textbf{0.00} & 48.51 &  & \textbf{13.42} & 1534.58 & 138.77 & 454.28 & 184.97 & 7200.00 \\
75 & 8 & 0.1 & 0.5 &  & \textbf{0.00} & 7.49 & 0.82 & 17.41 & 15.09 & 78.53 &  & \textbf{303.41} & 5897.76 & 3417.46 & 7200.00 & 7200.00 & 7200.00 \\
75 & 8 & 0.1 & 1.0 &  & 19.94 & 17.06 & 12.82 & 12.91 & 12.92 & \textbf{4.28} &  & 7200.00 & 7200.00 & 7200.00 & 7200.00 & 7200.00 & \textbf{5913.22} \\
75 & 8 & 0.3 & 0.1 &  & \textbf{0.00} & \textbf{0.00} & \textbf{0.00} & \textbf{0.00} & \textbf{0.00} & 48.73 &  & \textbf{12.22} & 1639.33 & 156.50 & 203.78 & 192.79 & 7200.00 \\
75 & 8 & 0.3 & 0.5 &  & \textbf{4.47} & 24.38 & 11.24 & 30.29 & 29.17 & 82.37 &  & \textbf{5898.63} & 7200.00 & 7200.00 & 7200.00 & 7200.00 & 7200.00 \\
75 & 8 & 0.3 & 1.0 &  & 24.83 & 29.90 & 25.21 & 25.66 & 24.66 & \textbf{17.89} &  & 7200.00 & 7200.00 & 7200.00 & 7200.00 & 7200.00 & 7200.00 \\
\rowcolor[HTML]{EFEFEF} 
75 & 12 & 0.1 & 0.1 &  & \textbf{0.00} & 43.47 & \textbf{0.00} & 8.60 & 5.73 & 60.12 &  & 1151.54 & 7200.00 & \textbf{486.97} & 5791.59 & 4397.36 & 7200.00 \\
\rowcolor[HTML]{EFEFEF} 
75 & 12 & 0.1 & 0.5 &  & 18.61 & 38.51 & \textbf{14.04} & 27.07 & 25.68 & 83.61 &  & 7200.00 & 7200.00 & 7200.00 & 7200.00 & 7200.00 & 7200.00 \\
\rowcolor[HTML]{EFEFEF} 
75 & 12 & 0.1 & 1.0 &  & 35.38 & 36.36 & 32.03 & 33.37 & 31.28 & \textbf{17.09} &  & 7200.00 & 7200.00 & 7200.00 & 7200.00 & 7200.00 & 7200.00 \\
\rowcolor[HTML]{EFEFEF} 
75 & 12 & 0.3 & 0.1 &  & \textbf{0.00} & 32.53 & \textbf{0.00} & 2.01 & \textbf{0.00} & 56.73 &  & 699.75 & 7200.00 & \textbf{287.67} & 2048.48 & 682.45 & 7200.00 \\
\rowcolor[HTML]{EFEFEF} 
75 & 12 & 0.3 & 0.5 &  & 31.92 & 38.87 & \textbf{26.73} & 37.94 & 36.79 & 87.84 &  & 7200.00 & 7200.00 & 7200.00 & 7200.00 & 7200.00 & 7200.00 \\
\rowcolor[HTML]{EFEFEF} 
75 & 12 & 0.3 & 1.0 &  & 36.74 & 46.50 & 33.06 & 33.71 & 32.09 & \textbf{22.08} &  & 7200.00 & 7200.00 & 7200.00 & 7200.00 & 7200.00 & 7200.00 \\
100 & 4 & 0.1 & 0.1 &  & \textbf{0.00} & \textbf{0.00} & \textbf{0.00} & \textbf{0.00} & \textbf{0.00} & 1.09 &  & \textbf{0.64} & 5.99 & 177.49 & 36.13 & 19.69 & 1663.17 \\
100 & 4 & 0.1 & 0.5 &  & \textbf{0.00} & \textbf{0.00} & \textbf{0.00} & \textbf{0.00} & \textbf{0.00} & 77.23 &  & \textbf{1.47} & 30.85 & 314.89 & 3098.62 & 2132.05 & 7200.00 \\
100 & 4 & 0.1 & 1.0 &  & \textbf{0.00} & \textbf{0.00} & \textbf{0.00} & \textbf{0.00} & \textbf{0.00} & 4.67 &  & 500.51 & \textbf{6.46} & 431.81 & 1237.74 & 673.01 & 5836.16 \\
100 & 4 & 0.3 & 0.1 &  & \textbf{0.00} & \textbf{0.00} & \textbf{0.00} & \textbf{0.00} & \textbf{0.00} & 5.84 &  & \textbf{0.36} & 3.39 & 69.54 & 16.54 & 11.78 & 1579.84 \\
100 & 4 & 0.3 & 0.5 &  & \textbf{0.00} & \textbf{0.00} & \textbf{0.00} & 2.97 & 1.92 & 78.02 &  & \textbf{2.47} & 91.93 & 318.28 & 4297.90 & 3750.83 & 7200.00 \\
100 & 4 & 0.3 & 1.0 &  & 0.21 & \textbf{0.00} & \textbf{0.00} & \textbf{0.00} & \textbf{0.00} & 6.16 &  & 1679.05 & \textbf{36.07} & 1188.05 & 2410.31 & 2957.87 & 7200.00 \\
\rowcolor[HTML]{EFEFEF} 
100 & 8 & 0.1 & 0.1 &  & \textbf{0.00} & 2.22 & \textbf{0.00} & 2.41 & \textbf{0.00} & 60.57 &  & \textbf{24.37} & 3328.77 & 184.96 & 2431.60 & 772.72 & 7200.00 \\
\rowcolor[HTML]{EFEFEF} 
100 & 8 & 0.1 & 0.5 &  & \textbf{0.90} & 13.43 & 7.06 & 19.17 & 17.93 & 85.26 &  & \textbf{3408.82} & 7200.00 & 7200.00 & 7200.00 & 7200.00 & 7200.00 \\
\rowcolor[HTML]{EFEFEF} 
100 & 8 & 0.1 & 1.0 &  & 25.08 & 32.20 & 22.55 & 23.08 & 21.33 & \textbf{2.16} &  & 7200.00 & 7200.00 & 7200.00 & 7200.00 & 7200.00 & \textbf{6561.62} \\
\rowcolor[HTML]{EFEFEF} 
100 & 8 & 0.3 & 0.1 &  & \textbf{0.00} & 0.66 & \textbf{0.00} & \textbf{0.00} & \textbf{0.00} & 60.95 &  & \textbf{26.80} & 3112.48 & 142.64 & 463.14 & 141.22 & 7200.00 \\
\rowcolor[HTML]{EFEFEF} 
100 & 8 & 0.3 & 0.5 &  & \textbf{9.40} & 24.04 & 13.86 & 26.25 & 24.77 & 87.96 &  & \textbf{5955.77} & 7200.00 & 7200.00 & 7200.00 & 7200.00 & 7200.00 \\
\rowcolor[HTML]{EFEFEF} 
100 & 8 & 0.3 & 1.0 &  & 23.24 & 25.55 & 22.07 & 22.89 & 21.88 & \textbf{13.64} &  & 7200.00 & 7200.00 & 7200.00 & 7200.00 & 7200.00 & \textbf{6099.95} \\
100 & 12 & 0.1 & 0.1 &  & 2.84 & 47.99 & \textbf{0.00} & 20.22 & 15.64 & 72.12 &  & 5774.24 & 7200.00 & \textbf{1562.89} & 7200.00 & 7053.93 & 7200.00 \\
100 & 12 & 0.1 & 0.5 &  & 23.42 & 38.72 & \textbf{19.86} & 28.80 & 27.96 & 89.31 &  & 7200.00 & 7200.00 & 7200.00 & 7200.00 & 7200.00 & 7200.00 \\
100 & 12 & 0.1 & 1.0 &  & 37.29 & 37.55 & 35.03 & 35.88 & 34.09 & \textbf{19.46} &  & 7200.00 & 7200.00 & 7200.00 & 7200.00 & 7200.00 & 7200.00 \\
100 & 12 & 0.3 & 0.1 &  & 3.57 & 33.99 & \textbf{0.00} & 6.48 & 6.85 & 71.39 &  & 4964.45 & 7200.00 & \textbf{1347.38} & 5057.35 & 5210.32 & 7200.00 \\
100 & 12 & 0.3 & 0.5 &  & 35.49 & 40.82 & \textbf{32.15} & 39.78 & 38.61 & 92.01 &  & 7200.00 & 7200.00 & 7200.00 & 7200.00 & 7200.00 & 7200.00 \\
100 & 12 & 0.3 & 1.0 &  & 32.92 & 72.74 & 30.98 & 31.71 & 30.06 & \textbf{23.89} &  & 7200.00 & 7200.00 & 7200.00 & 7200.00 & 7200.00 & 7200.00 \\
\rowcolor[HTML]{EFEFEF} 
100 & 16 & 0.1 & 0.1 &  & 43.18 & 60.62 & \textbf{11.27} & 33.45 & 29.34 & 77.93 &  & 7200.00 & 7200.00 & \textbf{5878.11} & 7200.00 & 7200.00 & 7200.00 \\
\rowcolor[HTML]{EFEFEF} 
100 & 16 & 0.1 & 0.5 &  & 35.41 & 58.47 & \textbf{30.08} & 34.33 & 33.09 & 92.28 &  & 7200.00 & 7200.00 & 7200.00 & 7200.00 & 7200.00 & 7200.00 \\
\rowcolor[HTML]{EFEFEF} 
100 & 16 & 0.1 & 1.0 &  & 33.90 & 100.00 & 31.71 & 33.14 & 31.14 & \textbf{20.17} &  & 7200.00 & 7200.00 & 7200.00 & 7200.00 & 7200.00 & 7200.00 \\
\rowcolor[HTML]{EFEFEF} 
100 & 16 & 0.3 & 0.1 &  & 64.02 & 71.55 & \textbf{7.35} & 29.36 & 27.62 & 78.97 &  & 7200.00 & 7200.00 & \textbf{5716.30} & 7200.00 & 7200.00 & 7200.00 \\
\rowcolor[HTML]{EFEFEF} 
100 & 16 & 0.3 & 0.5 &  & 64.27 & 66.41 & \textbf{46.52} & 50.86 & 50.39 & 93.69 &  & 7200.00 & 7200.00 & 7200.00 & 7200.00 & 7200.00 & 7200.00 \\
\rowcolor[HTML]{EFEFEF} 
100 & 16 & 0.3 & 1.0 &  & 58.98 & 86.10 & 43.57 & 44.49 & 42.89 & \textbf{21.14} &  & 7200.00 & 7200.00 & 7200.00 & 7200.00 & 7200.00 & 7200.00 \\ \bottomrule
\end{tabular}
}
\caption{Detailed results of the formulations $\bm{\mathrm{C}}$, $\bm{\mathrm{C}}^{+}_{e}$, $\bm{\mathrm{ICC}}$, $\bm{\mathrm{ICC+}}$, $\bm{\mathrm{LICC+}}$ and $\bm{\mathrm{CF}}$ for $\bm{\Gamma=0.9}$}
\label{tab:results_avg_0_9}
\end{table}

\newpage
\begin{table}[H]
\centering
\scalebox{0.585}{
\begin{tabular}{rrrrrrrrrrrrrrrrrrrr}
\toprule
\multicolumn{1}{c}{} & \multicolumn{1}{c}{} & \multicolumn{1}{c}{} & \multicolumn{1}{c}{} & \multicolumn{1}{c}{\textbf{}} & \multicolumn{7}{c}{\textbf{Avg. gap (\%)}} & \multicolumn{1}{c}{\textbf{}} & \multicolumn{7}{c}{\textbf{Avg. total time (s)}} \\ \cmidrule(lr){6-12} \cmidrule(l){14-20}
\multicolumn{1}{c}{\textbf{n}} & \multicolumn{1}{c}{\textbf{K}} & \multicolumn{1}{c}{$\boldsymbol{\beta}$} & \multicolumn{1}{c}{$\boldsymbol{\alpha}$} & \multicolumn{1}{c}{} & \multicolumn{1}{c}{$\bm{\mathrm{A}}$} & \multicolumn{1}{c}{$\bm{\mathrm{C}}$} & \multicolumn{1}{c}{$\bm{\mathrm{C}}^{+}_{e}$} & \multicolumn{1}{c}{$\bm{\mathrm{ICC}}$} & \multicolumn{1}{c}{$\bm{\mathrm{ICC+}}$} & \multicolumn{1}{c}{$\bm{\mathrm{LICC+}}$ } & \multicolumn{1}{c}{$\bm{\mathrm{CF}}$} & \multicolumn{1}{c}{} &  \multicolumn{1}{c}{$\bm{\mathrm{A}}$} & \multicolumn{1}{c}{$\bm{\mathrm{C}}$} & \multicolumn{1}{c}{$\bm{\mathrm{C}}^{+}_{e}$} & \multicolumn{1}{c}{$\bm{\mathrm{ICC}}$} & \multicolumn{1}{c}{$\bm{\mathrm{ICC+}}$} & \multicolumn{1}{c}{$\bm{\mathrm{LICC+}}$} & \multicolumn{1}{c}{$\bm{\mathrm{CF}}$} \\ \midrule
\rowcolor[HTML]{EFEFEF} 
50 & 4 & 0.1 & 0.1 &  & 10.0 & \textbf{0.0} & \textbf{0.0} & \textbf{0.0} & \textbf{0.0} & \textbf{0.0} & \textbf{0.0} &  & 1443.8 & \textbf{2.2} & 3.0 & 59.8 & 17.6 & 13.1 & 5.9 \\
\rowcolor[HTML]{EFEFEF} 
50 & 4 & 0.1 & 0.5 &  & 10.0 & \textbf{0.0} & \textbf{0.0} & \textbf{0.0} & \textbf{0.0} & \textbf{0.0} & \textbf{0.0} &  & 1831.8 & 70.2 & \textbf{4.7} & 244.9 & 342.7 & 188.7 & 12.9 \\
\rowcolor[HTML]{EFEFEF} 
50 & 4 & 0.1 & 1.0 &  & 17.3 & 0.0 & 0.0 & 0.0 & 0.0 & \textbf{0.0} & 0.0 &  & 2973.3 & 98.5 & \textbf{1.9} & 251.5 & 519.7 & 150.4 & 11.3 \\
\rowcolor[HTML]{EFEFEF} 
50 & 4 & 0.3 & 0.1 &  & \textbf{0.0} & \textbf{0.0} & \textbf{0.0} & \textbf{0.0} & \textbf{0.0} & \textbf{0.0} & \textbf{0.0} &  & 1.0 & \textbf{0.4} & 0.7 & 59.6 & 16.0 & 13.7 & 5.8 \\
\rowcolor[HTML]{EFEFEF} 
50 & 4 & 0.3 & 0.5 &  & \textbf{0.0} & \textbf{0.0} & \textbf{0.0} & \textbf{0.0} & 29.3 & 9.7 & \textbf{0.0} &  & 1087.5 & 1039.9 & \textbf{273.7} & 565.7 & 5819.4 & 2994.4 & 302.3 \\
\rowcolor[HTML]{EFEFEF} 
50 & 4 & 0.3 & 1.0 &  & 5.8 & \textbf{0.0} & \textbf{0.0} & \textbf{0.0} & 5.5 & \textbf{0.0} & \textbf{0.0} &  & 2585.3 & 944.8 & \textbf{29.1} & 711.8 & 3182.5 & 1200.8 & 52.8 \\
50 & 8 & 0.1 & 0.1 &  & 83.6 & 98.2 & 25.9 & 26.9 & \textbf{0.0} & \textbf{0.0} & \textbf{0.0} &  & 7200.0 & 7200.0 & 4928.1 & 6482.1 & 634.8 & 760.0 & \textbf{36.4} \\
50 & 8 & 0.1 & 0.5 &  & 89.6 & 99.2 & 41.8 & 54.0 & 44.0 & 44.1 & \textbf{0.0} &  & 7200.0 & 7200.0 & 5871.7 & 7200.0 & 6217.6 & 6276.6 & \textbf{1977.6} \\
50 & 8 & 0.1 & 1.0 &  & 90.6 & 97.8 & 25.7 & 48.8 & 40.0 & 32.3 & \textbf{0.0} &  & 7200.0 & 7200.0 & 5893.7 & 7200.0 & 7200.0 & 6558.8 & \textbf{728.6} \\
50 & 8 & 0.3 & 0.1 &  & 63.0 & 89.5 & 39.5 & 20.3 & \textbf{0.0} & \textbf{0.0} & \textbf{0.0} &  & 7200.0 & 7200.0 & 7122.8 & 4844.1 & 279.5 & 461.9 & \textbf{40.4} \\
50 & 8 & 0.3 & 0.5 &  & 90.8 & 96.3 & 66.2 & 75.0 & 71.9 & 71.8 & \textbf{52.2} &  & 7200.0 & 7200.0 & 7200.0 & 7200.0 & 7200.0 & 7200.0 & 7200.0 \\
50 & 8 & 0.3 & 1.0 &  & 85.1 & 95.2 & 61.6 & 65.6 & 64.1 & 60.2 & \textbf{30.1} &  & 7200.0 & 7200.0 & 7200.0 & 7200.0 & 7200.0 & 7200.0 & \textbf{6747.3} \\
\rowcolor[HTML]{EFEFEF} 
75 & 4 & 0.1 & 0.1 &  & 10.0 & \textbf{0.0} & \textbf{0.0} & \textbf{0.0} & 6.2 & \textbf{0.0} & \textbf{0.0} &  & 1454.4 & 13.4 & \textbf{5.9} & 81.1 & 1459.6 & 17.8 & 7.1 \\
\rowcolor[HTML]{EFEFEF} 
75 & 4 & 0.1 & 0.5 &  & 20.0 & \textbf{0.0} & \textbf{0.0} & 4.4 & 9.8 & 4.0 & \textbf{0.0} &  & 4887.5 & 293.8 & 132.8 & 1649.1 & 2879.7 & 1878.8 & \textbf{107.5} \\
\rowcolor[HTML]{EFEFEF} 
75 & 4 & 0.1 & 1.0 &  & 21.6 & \textbf{0.0} & \textbf{0.0} & \textbf{0.0} & 3.7 & 11.2 & \textbf{0.0} &  & 5679.1 & 1350.2 & \textbf{4.8} & 308.3 & 2277.4 & 2934.8 & 15.0 \\
\rowcolor[HTML]{EFEFEF} 
75 & 4 & 0.3 & 0.1 &  & \textbf{0.0} & \textbf{0.0} & \textbf{0.0} & \textbf{0.0} & \textbf{0.0} & \textbf{0.0} & \textbf{0.0} &  & 1.7 & \textbf{1.1} & 10.8 & 75.5 & 22.0 & 23.3 & 6.8 \\
\rowcolor[HTML]{EFEFEF} 
75 & 4 & 0.3 & 0.5 &  & 8.2 & \textbf{0.0} & \textbf{0.0} & 10.0 & 37.6 & 30.7 & 7.1 &  & 2716.8 & 1158.4 & \textbf{1037.9} & 3825.4 & 6803.0 & 6576.6 & 2694.5 \\
\rowcolor[HTML]{EFEFEF} 
75 & 4 & 0.3 & 1.0 &  & 22.2 & 5.6 & \textbf{0.0} & 4.4 & 22.0 & 8.3 & \textbf{0.0} &  & 5988.2 & 1618.1 & \textbf{56.8} & 2093.2 & 7200.0 & 5440.9 & 57.6 \\
75 & 8 & 0.1 & 0.1 &  & 93.2 & 99.7 & 70.9 & 56.8 & 28.1 & 41.7 & \textbf{8.8} &  & 7200.0 & 7200.0 & 7200.0 & 7200.0 & 6120.6 & 7200.0 & \textbf{3872.1} \\
75 & 8 & 0.1 & 0.5 &  & 95.4 & 99.0 & 65.3 & 71.1 & 65.6 & 61.0 & \textbf{15.2} &  & 7200.0 & 7200.0 & 7200.0 & 7200.0 & 7200.0 & 7200.0 & \textbf{4596.8} \\
75 & 8 & 0.1 & 1.0 &  & 95.6 & 96.6 & 65.4 & 60.7 & 65.2 & 56.2 & \textbf{2.4} &  & 7200.0 & 7200.0 & 7200.0 & 7200.0 & 7200.0 & 7200.0 & \textbf{2611.3} \\
75 & 8 & 0.3 & 0.1 &  & 89.7 & 94.1 & 69.5 & 47.4 & 20.2 & 27.6 & \textbf{0.0} &  & 7200.0 & 7200.0 & 7200.0 & 7200.0 & 4927.4 & 5938.2 & \textbf{1826.5} \\
75 & 8 & 0.3 & 0.5 &  & 96.3 & 97.6 & 85.6 & 84.6 & 80.5 & 79.4 & \textbf{68.1} &  & 7200.0 & 7200.0 & 7200.0 & 7200.0 & 7200.0 & 7200.0 & 7200.0 \\
75 & 8 & 0.3 & 1.0 &  & 94.3 & 97.9 & 70.5 & 76.9 & 77.8 & 76.0 & \textbf{54.8} &  & 7200.0 & 7200.0 & 7200.0 & 7200.0 & 7200.0 & 7200.0 & 7200.0 \\
\rowcolor[HTML]{EFEFEF} 
75 & 12 & 0.1 & 0.1 &  & 100.0 & 100.0 & 84.4 & 76.7 & 63.2 & 70.2 & \textbf{55.5} &  & 7200.0 & 7200.0 & 7200.0 & 7200.0 & 7200.0 & 7200.0 & 7200.0 \\
\rowcolor[HTML]{EFEFEF} 
75 & 12 & 0.1 & 0.5 &  & 100.0 & 100.0 & 94.2 & 85.6 & 79.3 & 76.6 & \textbf{65.7} &  & 7200.0 & 7200.0 & 7200.0 & 7200.0 & 7200.0 & 7200.0 & 7200.0 \\
\rowcolor[HTML]{EFEFEF} 
75 & 12 & 0.1 & 1.0 &  & 100.0 & 100.0 & 93.9 & 74.8 & 82.9 & 77.8 & \textbf{61.7} &  & 7200.0 & 7200.0 & 7200.0 & 7200.0 & 7200.0 & 7200.0 & 7200.0 \\
\rowcolor[HTML]{EFEFEF} 
75 & 12 & 0.3 & 0.1 &  & 100.0 & 100.0 & 87.4 & 75.5 & 64.8 & 71.0 & \textbf{54.9} &  & 7200.0 & 7200.0 & 7200.0 & 7200.0 & 7200.0 & 7200.0 & 7200.0 \\
\rowcolor[HTML]{EFEFEF} 
75 & 12 & 0.3 & 0.5 &  & 100.0 & 100.0 & 97.5 & 92.0 & 88.0 & 87.0 & \textbf{81.5} &  & 7200.0 & 7200.0 & 7200.0 & 7200.0 & 7200.0 & 7200.0 & 7200.0 \\
\rowcolor[HTML]{EFEFEF} 
75 & 12 & 0.3 & 1.0 &  & 100.0 & 100.0 & 98.2 & 85.1 & 89.0 & 85.1 & \textbf{75.0} &  & 7200.0 & 7200.0 & 7200.0 & 7200.0 & 7200.0 & 7200.0 & 7200.0 \\
100 & 4 & 0.1 & 0.1 &  & \textbf{0.0} & \textbf{0.0} & \textbf{0.0} & \textbf{0.0} & \textbf{0.0} & \textbf{0.0} & \textbf{0.0} &  & 276.7 & \textbf{3.4} & 9.3 & 117.9 & 83.0 & 95.8 & 13.8 \\
100 & 4 & 0.1 & 0.5 &  & 44.6 & \textbf{0.0} & \textbf{0.0} & 14.3 & 8.5 & 5.3 & \textbf{0.0} &  & 7200.0 & 702.6 & \textbf{265.6} & 3724.0 & 3717.2 & 2639.1 & 513.1 \\
100 & 4 & 0.1 & 1.0 &  & 57.5 & 3.7 & \textbf{0.0} & 4.8 & 6.7 & 11.7 & \textbf{0.0} &  & 7200.0 & 3301.2 & \textbf{13.1} & 2046.5 & 3471.9 & 3827.2 & 32.7 \\
100 & 4 & 0.3 & 0.1 &  & \textbf{0.0} & \textbf{0.0} & \textbf{0.0} & \textbf{0.0} & \textbf{0.0} & \textbf{0.0} & \textbf{0.0} &  & 12.5 & \textbf{3.0} & 208.3 & 123.9 & 57.6 & 128.3 & 20.8 \\
100 & 4 & 0.3 & 0.5 &  & 26.3 & \textbf{4.4} & 5.1 & 26.5 & 55.9 & 43.2 & 24.7 &  & 5904.1 & \textbf{2061.9} & 3459.8 & 5933.4 & 7200.0 & 7200.0 & 4500.2 \\
100 & 4 & 0.3 & 1.0 &  & 50.1 & 8.9 & \textbf{0.0} & 28.8 & 37.7 & 35.1 & \textbf{0.0} &  & 7200.0 & 5159.8 & 228.9 & 7200.0 & 7200.0 & 7200.0 & \textbf{91.6} \\
\rowcolor[HTML]{EFEFEF} 
100 & 8 & 0.1 & 0.1 &  & 86.0 & 93.6 & 74.9 & 52.9 & 39.5 & 48.2 & \textbf{16.8} &  & 7200.0 & 7200.0 & 7200.0 & 7200.0 & 7200.0 & 7200.0 & \textbf{4620.5} \\
\rowcolor[HTML]{EFEFEF} 
100 & 8 & 0.1 & 0.5 &  & 93.8 & 98.1 & 80.9 & 75.5 & 73.0 & 64.7 & \textbf{33.8} &  & 7200.0 & 7200.0 & 7200.0 & 7200.0 & 7200.0 & 7200.0 & \textbf{5281.8} \\
\rowcolor[HTML]{EFEFEF} 
100 & 8 & 0.1 & 1.0 &  & 93.3 & 96.5 & 86.8 & 62.3 & 67.9 & 59.7 & \textbf{9.1} &  & 7200.0 & 7200.0 & 7200.0 & 7200.0 & 7200.0 & 7200.0 & \textbf{4115.7} \\
\rowcolor[HTML]{EFEFEF} 
100 & 8 & 0.3 & 0.1 &  & 97.1 & 100.0 & 81.2 & 66.4 & 47.3 & 61.3 & \textbf{42.7} &  & 7200.0 & 7200.0 & 7200.0 & 7200.0 & 6311.6 & 7200.0 & \textbf{5831.8} \\
\rowcolor[HTML]{EFEFEF} 
100 & 8 & 0.3 & 0.5 &  & 99.1 & 100.0 & 93.9 & 91.9 & 88.1 & 87.5 & \textbf{81.9} &  & 7200.0 & 7200.0 & 7200.0 & 7200.0 & 7200.0 & 7200.0 & 7200.0 \\
\rowcolor[HTML]{EFEFEF} 
100 & 8 & 0.3 & 1.0 &  & 98.5 & 98.8 & 89.5 & 87.9 & 89.6 & 85.4 & \textbf{67.4} &  & 7200.0 & 7200.0 & 7200.0 & 7200.0 & 7200.0 & 7200.0 & 7200.0 \\
100 & 12 & 0.1 & 0.1 &  & 100.0 & 100.0 & 87.7 & 82.8 & 76.4 & 76.4 & \textbf{68.4} &  & 7200.0 & 7200.0 & 7200.0 & 7200.0 & 7200.0 & 7200.0 & 7200.0 \\
100 & 12 & 0.1 & 0.5 &  & 100.0 & 100.0 & 94.6 & 87.3 & 81.6 & 76.5 & \textbf{67.8} &  & 7200.0 & 7200.0 & 7200.0 & 7200.0 & 7200.0 & 7200.0 & 7200.0 \\
100 & 12 & 0.1 & 1.0 &  & 100.0 & 100.0 & 94.2 & 79.4 & 87.5 & 79.6 & \textbf{66.7} &  & 7200.0 & 7200.0 & 7200.0 & 7200.0 & 7200.0 & 7200.0 & 7200.0 \\
100 & 12 & 0.3 & 0.1 &  & 100.0 & 100.0 & 100.0 & 82.9 & 77.1 & 77.8 & \textbf{69.4} &  & 7200.0 & 7200.0 & 7200.0 & 7200.0 & 7200.0 & 7200.0 & 7200.0 \\
100 & 12 & 0.3 & 0.5 &  & 100.0 & 100.0 & 100.0 & 95.2 & 92.0 & 91.0 & \textbf{88.9} &  & 7200.0 & 7200.0 & 7200.0 & 7200.0 & 7200.0 & 7200.0 & 7200.0 \\
100 & 12 & 0.3 & 1.0 &  & 100.0 & 100.0 & 100.0 & 90.3 & 95.1 & 90.9 & \textbf{85.3} &  & 7200.0 & 7200.0 & 7200.0 & 7200.0 & 7200.0 & 7200.0 & 7200.0 \\
\rowcolor[HTML]{EFEFEF} 
100 & 16 & 0.1 & 0.1 &  & 100.0 & 100.0 & 100.0 & 93.0 & 80.3 & 80.3 & \textbf{73.6} &  & 7200.0 & 7200.0 & 7200.0 & 7200.0 & 7200.0 & 7200.0 & 7200.0 \\
\rowcolor[HTML]{EFEFEF} 
100 & 16 & 0.1 & 0.5 &  & 100.0 & 100.0 & 100.0 & 93.2 & 88.9 & 86.6 & \textbf{82.5} &  & 7200.0 & 7200.0 & 7200.0 & 7200.0 & 7200.0 & 7200.0 & 7200.0 \\
\rowcolor[HTML]{EFEFEF} 
100 & 16 & 0.1 & 1.0 &  & 100.0 & 100.0 & 100.0 & 86.8 & 98.7 & 86.4 & \textbf{81.1} &  & 7200.0 & 7200.0 & 7200.0 & 7200.0 & 7200.0 & 7200.0 & 7200.0 \\
\rowcolor[HTML]{EFEFEF} 
100 & 16 & 0.3 & 0.1 &  & 100.0 & 100.0 & 100.0 & 86.9 & 81.8 & 81.7 & \textbf{74.8} &  & 7200.0 & 7200.0 & 7200.0 & 7200.0 & 7200.0 & 7200.0 & 7200.0 \\
\rowcolor[HTML]{EFEFEF} 
100 & 16 & 0.3 & 0.5 &  & 100.0 & 100.0 & 100.0 & 96.1 & 93.2 & 92.3 & \textbf{90.8} &  & 7200.0 & 7200.0 & 7200.0 & 7200.0 & 7200.0 & 7200.0 & 7200.0 \\
\rowcolor[HTML]{EFEFEF} 
100 & 16 & 0.3 & 1.0 &  & 100.0 & 100.0 & 100.0 & 91.0 & 99.2 & 91.6 & \textbf{88.1} &  & 7200.0 & 7200.0 & 7200.0 & 7200.0 & 7200.0 & 7200.0 & 7200.0 \\ \bottomrule
\end{tabular}
}
\caption{Detailed results of the formulations $\bm{\mathrm{A}}$, $\bm{\mathrm{C}}$, $\bm{\mathrm{C}}^{+}_{e}$, $\bm{\mathrm{ICC}}$, $\bm{\mathrm{ICC+}}$, $\bm{\mathrm{LICC+}}$ and $\bm{\mathrm{CF}}$ for $\boldsymbol{\Gamma=1.0}$}
\label{tab:results_avg_1_0}
\end{table}

\newpage
\begin{table}[H]
\centering
\scalebox{0.585}{
\begin{tabular}{rrrrrrrrrrrrrrrrrr}
\toprule
\multicolumn{1}{c}{} & \multicolumn{1}{c}{} & \multicolumn{1}{c}{} & \multicolumn{1}{c}{} & \multicolumn{1}{c}{\textbf{}} & \multicolumn{6}{c}{\textbf{Avg. gap (\%)}} & \multicolumn{1}{c}{\textbf{}} & \multicolumn{6}{c}{\textbf{Avg. total time (s)}} \\ \cmidrule(lr){6-11} \cmidrule(l){13-18} 
\multicolumn{1}{c}{\textbf{n}} & \multicolumn{1}{c}{\textbf{K}} & \multicolumn{1}{c}{$\boldsymbol{\beta}$} & \multicolumn{1}{c}{$\boldsymbol{\alpha}$} & \multicolumn{1}{c}{} & \multicolumn{1}{c}{$\bm{\mathrm{C}}$} & \multicolumn{1}{c}{$\bm{\mathrm{C}}^{+}_{e}$} & \multicolumn{1}{c}{$\bm{\mathrm{ICC}}$} & \multicolumn{1}{c}{$\bm{\mathrm{ICC+}}$} & \multicolumn{1}{c}{$\bm{\mathrm{LICC+}}$ } & \multicolumn{1}{c}{$\bm{\mathrm{CF}}$} & \multicolumn{1}{c}{} & \multicolumn{1}{c}{$\bm{\mathrm{C}}$} & \multicolumn{1}{c}{$\bm{\mathrm{C}}^{+}_{e}$} & \multicolumn{1}{c}{$\bm{\mathrm{ICC}}$} & \multicolumn{1}{c}{$\bm{\mathrm{ICC+}}$} & \multicolumn{1}{c}{$\bm{\mathrm{LICC+}}$} & \multicolumn{1}{c}{$\bm{\mathrm{CF}}$} \\ \midrule
\rowcolor[HTML]{EFEFEF} 
50 & 4 & 0.1 & 0.1 &  & 20.00 & \textbf{0.00} & \textbf{0.00} & \textbf{0.00} & \textbf{0.00} & \textbf{0.00} &  & 1507.00 & \textbf{0.41} & 51.54 & 14.67 & 156.03 & 5.11 \\
\rowcolor[HTML]{EFEFEF} 
50 & 4 & 0.1 & 0.5 &  & 8.33 & \textbf{0.00} & \textbf{0.00} & \textbf{0.00} & \textbf{0.00} & \textbf{0.00} &  & 1729.82 & \textbf{1.58} & 171.53 & 22.58 & 18.18 & 5.36 \\
\rowcolor[HTML]{EFEFEF} 
50 & 4 & 0.1 & 1.0 &  & 8.33 & \textbf{0.00} & \textbf{0.00} & \textbf{0.00} & \textbf{0.00} & \textbf{0.00} &  & 1496.90 & \textbf{0.29} & 14.18 & 31.64 & 10.57 & 5.29 \\
\rowcolor[HTML]{EFEFEF} 
50 & 4 & 0.3 & 0.1 &  & \textbf{0.00} & \textbf{0.00} & \textbf{0.00} & 2.00 & \textbf{0.00} & \textbf{0.00} &  & 122.75 & \textbf{0.33} & 75.73 & 1798.93 & 1378.67 & 5.01 \\
\rowcolor[HTML]{EFEFEF} 
50 & 4 & 0.3 & 0.5 &  & 10.00 & \textbf{0.00} & 1.33 & 14.87 & \textbf{0.00} & \textbf{0.00} &  & 3547.04 & \textbf{3.51} & 1604.75 & 2962.73 & 72.93 & 6.14 \\
\rowcolor[HTML]{EFEFEF} 
50 & 4 & 0.3 & 1.0 &  & 13.00 & \textbf{0.00} & \textbf{0.00} & \textbf{0.00} & 11.25 & \textbf{0.00} &  & 3360.42 & \textbf{4.91} & 235.82 & 118.29 & 2896.49 & 6.28 \\
50 & 8 & 0.1 & 0.1 &  & 100.00 & 20.77 & 54.30 & 0.97 & 1.21 & \textbf{0.00} &  & 7200.00 & 4529.79 & 7200.00 & 2911.30 & 3093.84 & \textbf{11.28} \\
50 & 8 & 0.1 & 0.5 &  & 100.00 & 25.82 & 55.41 & 19.78 & 8.89 & \textbf{0.00} &  & 7200.00 & 5235.15 & 7200.00 & 4179.85 & 4560.54 & \textbf{35.03} \\
50 & 8 & 0.1 & 1.0 &  & 100.00 & 30.49 & 26.30 & \textbf{0.00} & 12.35 & \textbf{0.00} &  & 7200.00 & 4649.12 & 7200.00 & 2119.95 & 4712.71 & \textbf{64.84} \\
50 & 8 & 0.3 & 0.1 &  & 100.00 & 33.69 & 51.56 & 2.53 & 1.13 & \textbf{0.00} &  & 7200.00 & 5468.89 & 7200.00 & 5803.24 & 3193.92 & \textbf{14.17} \\
50 & 8 & 0.3 & 0.5 &  & 100.00 & 61.57 & 70.30 & 59.67 & 40.92 & \textbf{0.00} &  & 7200.00 & 7200.00 & 7200.00 & 7200.00 & 6403.32 & \textbf{288.39} \\
50 & 8 & 0.3 & 1.0 &  & 100.00 & 72.53 & 53.08 & 31.11 & 36.14 & \textbf{0.00} &  & 7200.00 & 7200.00 & 7200.00 & 6521.26 & 6530.25 & \textbf{209.56} \\
\rowcolor[HTML]{EFEFEF} 
75 & 4 & 0.1 & 0.1 &  & 40.00 & \textbf{0.00} & \textbf{0.00} & \textbf{0.00} & \textbf{0.00} & \textbf{0.00} &  & 4289.03 & \textbf{0.11} & 110.66 & 22.08 & 12.95 & 5.28 \\
\rowcolor[HTML]{EFEFEF} 
75 & 4 & 0.1 & 0.5 &  & 60.00 & \textbf{0.00} & 2.22 & 2.50 & \textbf{0.00} & \textbf{0.00} &  & 5349.42 & \textbf{0.78} & 1579.96 & 1465.81 & 26.09 & 5.06 \\
\rowcolor[HTML]{EFEFEF} 
75 & 4 & 0.1 & 1.0 &  & \textbf{0.00} & \textbf{0.00} & \textbf{0.00} & \textbf{0.00} & \textbf{0.00} & \textbf{0.00} &  & 1995.68 & \textbf{0.43} & 29.10 & 32.66 & 167.76 & 5.05 \\
\rowcolor[HTML]{EFEFEF} 
75 & 4 & 0.3 & 0.1 &  & 20.00 & \textbf{0.00} & \textbf{0.00} & 1.18 & \textbf{0.00} & \textbf{0.00} &  & 1549.71 & 34.14 & 121.69 & 1454.83 & 17.42 & \textbf{5.33} \\
\rowcolor[HTML]{EFEFEF} 
75 & 4 & 0.3 & 0.5 &  & 32.22 & \textbf{0.00} & 12.86 & 5.60 & \textbf{0.00} & \textbf{0.00} &  & 4486.56 & 487.90 & 2964.17 & 2118.72 & 565.11 & \textbf{14.71} \\
\rowcolor[HTML]{EFEFEF} 
75 & 4 & 0.3 & 1.0 &  & 30.59 & \textbf{0.00} & 8.87 & 13.08 & \textbf{0.00} & \textbf{0.00} &  & 3291.26 & 87.60 & 2935.24 & 1559.31 & 252.76 & \textbf{13.86} \\
75 & 8 & 0.1 & 0.1 &  & 100.00 & 53.97 & 68.72 & 5.24 & 5.11 & \textbf{0.00} &  & 7200.00 & 7200.00 & 7200.00 & 3972.00 & 3130.23 & \textbf{50.11} \\
75 & 8 & 0.1 & 0.5 &  & 100.00 & 48.71 & 68.52 & 41.45 & 18.64 & \textbf{0.00} &  & 7200.00 & 7200.00 & 7200.00 & 5895.17 & 6370.93 & \textbf{72.29} \\
75 & 8 & 0.1 & 1.0 &  & 100.00 & 54.59 & 43.45 & 17.05 & 26.47 & \textbf{0.00} &  & 7200.00 & 7200.00 & 7200.00 & 4496.17 & 5958.33 & \textbf{73.38} \\
75 & 8 & 0.3 & 0.1 &  & 100.00 & 59.10 & 70.76 & 5.72 & 8.78 & \textbf{0.00} &  & 7200.00 & 7200.00 & 7200.00 & 5097.19 & 6172.99 & \textbf{125.00} \\
75 & 8 & 0.3 & 0.5 &  & 100.00 & 77.76 & 75.71 & 60.87 & 43.59 & \textbf{0.00} &  & 7200.00 & 7200.00 & 7200.00 & 7200.00 & 6181.48 & \textbf{221.21} \\
75 & 8 & 0.3 & 1.0 &  & 100.00 & 76.13 & 51.55 & 46.12 & 45.38 & \textbf{0.00} &  & 7200.00 & 7200.00 & 7200.00 & 6340.83 & 7200.00 & \textbf{835.33} \\
\rowcolor[HTML]{EFEFEF} 
75 & 12 & 0.1 & 0.1 &  & 100.00 & 79.29 & 85.29 & 68.55 & 67.83 & \textbf{40.73} &  & 7200.00 & 7200.00 & 7200.00 & 7200.00 & 7200.00 & 7200.00 \\
\rowcolor[HTML]{EFEFEF} 
75 & 12 & 0.1 & 0.5 &  & 100.00 & 85.22 & 87.57 & 77.66 & 75.87 & \textbf{28.84} &  & 7200.00 & 7200.00 & 7200.00 & 7200.00 & 7200.00 & \textbf{6502.52} \\
\rowcolor[HTML]{EFEFEF} 
75 & 12 & 0.1 & 1.0 &  & 100.00 & 82.22 & 73.11 & 73.20 & 73.20 & \textbf{28.44} &  & 7200.00 & 7200.00 & 7200.00 & 7200.00 & 7200.00 & 7200.00 \\
\rowcolor[HTML]{EFEFEF} 
75 & 12 & 0.3 & 0.1 &  & 100.00 & 81.11 & 86.17 & 69.00 & 72.28 & \textbf{49.83} &  & 7200.00 & 7200.00 & 7200.00 & 7200.00 & 7200.00 & 7200.00 \\
\rowcolor[HTML]{EFEFEF} 
75 & 12 & 0.3 & 0.5 &  & 100.00 & 87.82 & 90.29 & 81.19 & 80.30 & \textbf{60.25} &  & 7200.00 & 7200.00 & 7200.00 & 7200.00 & 7200.00 & 7200.00 \\
\rowcolor[HTML]{EFEFEF} 
75 & 12 & 0.3 & 1.0 &  & 100.00 & 89.59 & 80.70 & 83.45 & 82.62 & \textbf{62.68} &  & 7200.00 & 7200.00 & 7200.00 & 7200.00 & 7200.00 & 7200.00 \\
100 & 4 & 0.1 & 0.1 &  & 80.00 & \textbf{0.00} & \textbf{0.00} & \textbf{0.00} & \textbf{0.00} & \textbf{0.00} &  & 5760.21 & \textbf{0.70} & 132.33 & 30.62 & 17.84 & 5.15 \\
100 & 4 & 0.1 & 0.5 &  & 80.00 & \textbf{0.00} & 20.67 & \textbf{0.00} & \textbf{0.00} & \textbf{0.00} &  & 5772.08 & 7.40 & 3030.47 & 299.81 & 140.51 & \textbf{5.77} \\
100 & 4 & 0.1 & 1.0 &  & 80.00 & \textbf{0.00} & \textbf{0.00} & \textbf{0.00} & \textbf{0.00} & \textbf{0.00} &  & 5761.34 & \textbf{5.57} & 136.07 & 208.41 & 55.43 & 6.03 \\
100 & 4 & 0.3 & 0.1 &  & 100.00 & \textbf{0.00} & 10.00 & 16.78 & 4.44 & \textbf{0.00} &  & 7200.00 & 224.34 & 2646.24 & 5864.04 & 2998.54 & \textbf{6.97} \\
100 & 4 & 0.3 & 0.5 &  & 100.00 & \textbf{0.00} & 52.37 & 17.31 & \textbf{0.00} & \textbf{0.00} &  & 7200.00 & 175.89 & 7200.00 & 5586.52 & 599.96 & \textbf{9.08} \\
100 & 4 & 0.3 & 1.0 &  & 86.65 & \textbf{0.00} & 7.15 & 0.83 & \textbf{0.00} & \textbf{0.00} &  & 7200.00 & 109.82 & 5785.19 & 1934.88 & 96.27 & \textbf{9.10} \\
\rowcolor[HTML]{EFEFEF} 
100 & 8 & 0.1 & 0.1 &  & 100.00 & 66.15 & 75.48 & 21.40 & 18.66 & \textbf{0.00} &  & 7200.00 & 7200.00 & 7200.00 & 5123.75 & 5865.88 & \textbf{157.61} \\
\rowcolor[HTML]{EFEFEF} 
100 & 8 & 0.1 & 0.5 &  & 100.00 & 47.84 & 73.78 & 49.98 & 26.69 & \textbf{0.00} &  & 7200.00 & 7200.00 & 7200.00 & 7200.00 & 4650.75 & \textbf{45.49} \\
\rowcolor[HTML]{EFEFEF} 
100 & 8 & 0.1 & 1.0 &  & 100.00 & 64.56 & 45.97 & 40.53 & 45.11 & \textbf{0.00} &  & 7200.00 & 7200.00 & 7200.00 & 6405.35 & 7200.00 & \textbf{93.72} \\
\rowcolor[HTML]{EFEFEF} 
100 & 8 & 0.3 & 0.1 &  & 100.00 & 79.84 & 79.02 & 53.07 & 50.86 & \textbf{0.00} &  & 7200.00 & 7200.00 & 7200.00 & 7200.00 & 7200.00 & \textbf{971.93} \\
\rowcolor[HTML]{EFEFEF} 
100 & 8 & 0.3 & 0.5 &  & 100.00 & 83.51 & 87.83 & 77.25 & 75.60 & \textbf{0.00} &  & 7200.00 & 7200.00 & 7200.00 & 7200.00 & 7200.00 & \textbf{2434.96} \\
\rowcolor[HTML]{EFEFEF} 
100 & 8 & 0.3 & 1.0 &  & 100.00 & 85.78 & 73.83 & 70.02 & 71.28 & \textbf{5.09} &  & 7200.00 & 7200.00 & 7200.00 & 7200.00 & 7200.00 & \textbf{2129.35} \\
100 & 12 & 0.1 & 0.1 &  & 100.00 & 82.32 & 89.05 & 72.40 & 72.71 & \textbf{56.95} &  & 7200.00 & 7200.00 & 7200.00 & 7200.00 & 7200.00 & 7200.00 \\
100 & 12 & 0.1 & 0.5 &  & 100.00 & 83.47 & 88.90 & 83.03 & 75.29 & \textbf{40.68} &  & 7200.00 & 7200.00 & 7200.00 & 7200.00 & 7200.00 & 7200.00 \\
100 & 12 & 0.1 & 1.0 &  & 100.00 & 85.19 & 72.05 & 75.15 & 74.19 & \textbf{36.60} &  & 7200.00 & 7200.00 & 7200.00 & 7200.00 & 7200.00 & 7200.00 \\
100 & 12 & 0.3 & 0.1 &  & 100.00 & 87.33 & 96.01 & 76.95 & 76.94 & \textbf{69.47} &  & 7200.00 & 7200.00 & 7200.00 & 7200.00 & 7200.00 & 7200.00 \\
100 & 12 & 0.3 & 0.5 &  & 100.00 & 93.26 & 94.04 & 88.51 & 87.27 & \textbf{71.57} &  & 7200.00 & 7200.00 & 7200.00 & 7200.00 & 7200.00 & 7200.00 \\
100 & 12 & 0.3 & 1.0 &  & 100.00 & 92.22 & 85.26 & 88.78 & 86.66 & \textbf{74.29} &  & 7200.00 & 7200.00 & 7200.00 & 7200.00 & 7200.00 & 7200.00 \\
\rowcolor[HTML]{EFEFEF} 
100 & 16 & 0.1 & 0.1 &  & 100.00 & 100.00 & 93.91 & 91.38 & 89.70 & \textbf{76.24} &  & 7200.00 & 7200.00 & 7200.00 & 7200.00 & 7200.00 & 7200.00 \\
\rowcolor[HTML]{EFEFEF} 
100 & 16 & 0.1 & 0.5 &  & 100.00 & 100.00 & 99.61 & 99.45 & 98.24 & \textbf{70.85} &  & 7200.00 & 7200.00 & 7200.00 & 7200.00 & 7200.00 & 7200.00 \\
\rowcolor[HTML]{EFEFEF} 
100 & 16 & 0.1 & 1.0 &  & 100.00 & 100.00 & 99.78 & 99.86 & 99.79 & \textbf{73.00} &  & 7200.00 & 7200.00 & 7200.00 & 7200.00 & 7200.00 & 7200.00 \\
\rowcolor[HTML]{EFEFEF} 
100 & 16 & 0.3 & 0.1 &  & 100.00 & 100.00 & 98.39 & 91.25 & 88.04 & \textbf{78.15} &  & 7200.00 & 7200.00 & 7200.00 & 7200.00 & 7200.00 & 7200.00 \\
\rowcolor[HTML]{EFEFEF} 
100 & 16 & 0.3 & 0.5 &  & 100.00 & 100.00 & 100.00 & 100.00 & 100.00 & \textbf{79.86} &  & 7200.00 & 7200.00 & 7200.00 & 7200.00 & 7200.00 & 7200.00 \\
\rowcolor[HTML]{EFEFEF} 
100 & 16 & 0.3 & 1.0 &  & 100.00 & 100.00 & 100.00 & 100.00 & 100.00 & \textbf{79.82} &  & 7200.00 & 7200.00 & 7200.00 & 7200.00 & 7200.00 & 7200.00 \\ \bottomrule
\end{tabular}
}
\caption{Detailed results of the formulations $\bm{\mathrm{C}}$, $\bm{\mathrm{C}}^{+}_{e}$, $\bm{\mathrm{ICC}}$, $\bm{\mathrm{ICC+}}$, $\bm{\mathrm{LICC+}}$ and $\bm{\mathrm{CF}}$ for $\boldsymbol{\Gamma=1.1}$}
\label{tab:results_avg_1_1}
\end{table}

\end{APPENDICES}

\newpage

\end{document}